\documentclass[12pt]{article}
\usepackage{amsmath}
\usepackage{latexsym}
\usepackage{amssymb}
\usepackage{a4wide}
\newtheorem{thm}{Theorem}[section]
\newtheorem{la}[thm]{Lemma}
\newtheorem{Defn}[thm]{Definition}
\newtheorem{Remark}[thm]{Remark}
\newtheorem{Note}[thm]{Note}
\newtheorem{prop}[thm]{Proposition}

\newtheorem{Example}[thm]{Example}
\newtheorem{Examples}[thm]{Examples}
\newtheorem{Problems}[thm]{Problems}

\newtheorem{Problem}[thm]{Problem}
\newtheorem{Number}[thm]{\!\!}
\newenvironment{defn}{\begin{Defn}\rm}{\end{Defn}}

\newenvironment{rem}{\begin{Remark}\rm}{\end{Remark}}
\newenvironment{numba}{\begin{Number}\rm}{\end{Number}}
\newenvironment{proof}{{\noindent\bf Proof.}}%
                  {\nopagebreak\hspace*{\fill}$\Box$\medskip\par}   
\newcommand{\Punkt}{\nopagebreak\hspace*{\fill}$\Box$}
\newcommand{\wb}{\overline}

\newcommand{\wt}{\widetilde}

\newcommand{\n}{\rm}

\newcommand{\impl}{\Rightarrow}

\newcommand{\mto}{\mapsto}

\newcommand{\ve}{\varepsilon}

\newcommand{\N}{{\mathbb N}}

\newcommand{\R}{{\mathbb R}}

\newcommand{\Q}{{\mathbb Q}}

\newcommand{\Z}{{\mathbb Z}}

\newcommand{\cO}{{\cal O}}

\newcommand{\one}{\mbox{\rm \bf 1}}

\newcommand{\sub}{\subseteq}

\newcommand{\im}{\mbox{\n im}}

\newcommand{\cL}{{\cal L}}

\newcommand{\Spann}{\mbox{\n span}}

\newcommand{\ev}{\mbox{\n ev}}

\newcommand{\sbull}{{\scriptscriptstyle \bullet}}

\newcommand{\absconv}{{\mbox{{\rm absconv}}}}

\newcommand{\half}{{\textstyle \frac{1}{2}}}
\newcommand{\quart}{{\textstyle \frac{1}{4}}}
\newcommand{\Lip}{\mbox{${\ell} {i} {p}$}}
\begin{document}
\begin{center}
{\Large \bf Conveniently H\"{o}lder Homomorphisms
are Smooth\\[2mm]
in the Convenient Sense}\vspace{4.5 mm}\\
{\bf Helge Gl\"{o}ckner}\vspace{.9mm}
\end{center}
\noindent{\bf Abstract.\/}
We show that
every ``conveniently H\"{o}lder''
homomorphism
between Lie groups
in the sense of convenient differential
calculus is smooth (in the convenient sense).
In particular,
every $\Lip^0$-homomorphism is smooth.\\[3mm]
{\footnotesize
{\bf AMS Subject Classification.}
22E65 (main),
26E15, 26E20, 46T20, 58C20.\\[3mm]
{\bf Keywords and Phrases.}
Infinite-dimensional Lie group, homomorphism,
Lipschitz condition, H\"{o}lder condition,
Taylor expansion, differentiability, smoothness,
convenient differential calculus.}
\section*{Introduction}
In a preprint from 1982,
John Milnor formulated various
absolutely fundamental
open problems concerning
infinite-dimensional
Lie groups~\cite{Mi1}.
Our investigations
are related to Milnor's
question:
{\em Is a continuous homomorphism
between Lie groups necessarily smooth\/}\,?,
which refers to smooth
Lie groups modelled on complete
locally convex spaces,
based on smooth maps
in the sense of Michal-Bastiani
(Keller's $C^\infty_c$-maps).
While the answer to this question
is still unknown,
some progress has been
made recently:
every {\em H\"{o}lder\/} continuous
homomorphism is smooth~\cite{HOL}.
In particular, every Lipschitz
continuous homomorphism
is smooth.
The goal of this article
is to establish analogous
results in the framework
of infinite-dimensional analysis
and Lie theory
known as Convenient Differential Calculus
(see \cite{FaK}, \cite{KaM}).
In this setting, a map
is called $\Lip^0$ (or $\cL ip^0$)
if it takes smooth curves to Lipschitz curves.
Instead of continuous maps,
one considers $\Lip^0$-maps
as the more adequate fundamental notion here,
because the $\Lip^0$-property
is a purely bornological concept.
Our main result is the following (Theorem~\ref{main}):\\[3mm]
{\bf Main Theorem.}
{\em Let $G$ and $H$ be Lie groups in the sense
of convenient differential calculus
and $f\!: G\to H$ be a homomorphism.
If $f$ is $\Lip^0$,
then $f$ is smooth in the convenient
sense. More generally,
this conclusion remains valid if
$f$ is conveniently H\"{o}lder.}\\[3mm]
Here, a map $f$
is called {\em conveniently
H\"{o}lder\/}
if it is $h_\alpha$
for some $\alpha\in \;]0,1]$
in the sense that $f$ takes smooth curves
to H\"{o}lder continuous curves of H\"{o}lder
exponent~$\alpha$ (thus $h_1=\Lip^0$).\\[3mm]
{\bf Strategy of proof.}
To establish
our Main Theorem,
the strategy is to
encode smoothness of homomorphisms
in a suitable differentiability property
at the identity,
which can be checked in the
conveniently H\"{o}lder case.
Our starting point
is the (trivial) observation
that every $\Lip^1$-homomorphism
$f\!: G\to H$ is
smooth in the convenient sense
(Lemma~\ref{C1smooth}).
The proof then proceeds
in two main steps:
\begin{itemize}
\item
First, we show that
a homomorphism
$f\!: G\to H$ is $\Lip^1$
(and hence smooth)
provided it is ``curve differentiable''
at~$1$.
\item
The difficult task, then,
is to show
that
every conveniently
H\"{o}lder homomorphism
is curve
differentiable at~$1$.
What we actually establish is
{\em bornological\/} curve differentiability
at~$1$, a stronger (but more tangible)
property.
\end{itemize}
We remark that smoothness of H\"{o}lder
continuous homomorphisms
in the setting of Keller's $C^\infty_c$-theory
is proved in~\cite{HOL}
in two analogous steps.
The appropriate notion of
differentiability at a point
used there is total differentiability.
In the framework
of convenient differential calculus,
differentiability at a point
has not been considered much
in the literature.
It was therefore necessary
to develop various new concepts.
Curve differentiability
and bornological curve differentiability,
which we introduce here,
serve us as efficient tools
for the discussion of homomorphisms.
Once the basic facts concerning
these differentiability properties
are established,
Step~1 (as described above)
is easily performed.
The proof of Step~2 is rather
technical and much more difficult.
However, one central idea of the proof
is easily explained on an informal level,
and we describe it now.
We recommend to keep this basic idea
in the back of one's mind when reading Section~\ref{hardlabor}.
To shorten formulas and
increase the readability,
let us identify an open identity neighbourhood
of $H$ with a $c^\infty$-open $0$-neighbourhood
$V\sub L(H)$ for the moment
(such that $0$ becomes the identity element
of~$H$). Likewise, we identify an open identity neighbourhood
in~$G$ with a $c^\infty$-open $0$-neighbourhood
$U\sub L(G)$, which we choose so small that $f(U)\sub V$.\\[3mm]
{\bf The core idea.}
To establish bornological
curve differentiability
of~$f$ at~$1$, (among other things)
we need to show that $(f\circ \gamma)'(0)$ exists,
for each smooth curve
$\gamma\!: \R\to U$
such that $\gamma(0)=0$.
We now explain how a candidate
for $(f\circ \gamma)'(0)$ can be obtained.
The idea is to exploit
the first
order Taylor expansion
$x^2y=2x+y+R(x,y)$
of the map $(x,y)\mto x^2y\in V\sub L(H)$
(defined on some $c^\infty$-open $(0,0)$-neighbourhood
in $V\times V$).
For sufficiently small $t$, we have
\begin{eqnarray*}
f(\gamma(t)) &= & f(\gamma(\half t))^2 \, f(\gamma(\half t)^{-2}\gamma(t))\\
& = & 2\, f(\gamma(\half t))\;+\;f(\gamma(\half t)^{-2}\gamma(t))
\;+\; R\big(f(\gamma(\half t)),\, f(\gamma(\half t)^{-2}\gamma(t))\big)
\end{eqnarray*}
and thus
$2 f(\gamma(\half t))=
f(\gamma(t))
-f(\gamma(\half t)^{-2}\gamma(t))
-R\big(f(\gamma(\half t)),f(\gamma(\half t)^{-2}\gamma(t))\big)$.
Hence
\begin{eqnarray*}
4 f(\gamma(\quart t)) & = &
2 f(\gamma(\half t))
-2 f(\gamma(\quart t)^{-2}\gamma(\half t))
-2 R\big(f(\gamma(\quart t)),f(\gamma(\quart t)^{-2}\gamma(\half
t))\big)\\
&=&
f(\gamma(t))
-f(\gamma(\half t)^{-2}\gamma(t))
-R\big(f(\gamma(\half t)),f(\gamma(\half t)^{-2}\gamma(t))\big)\\
& & \qquad\quad  \;\, -\,2 f(\gamma(\quart t)^{-2}\gamma(\half t))
-2 R\big(f(\gamma(\quart t)),f(\gamma(\quart t)^{-2}\gamma(\half
t))\big)\,,
\end{eqnarray*}
using the preceding formula twice.
Repeating this argument,
we obtain
\begin{eqnarray}
\hspace*{-2cm}\frac{f(\gamma(2^{-n}t))}{2^{-n}}
&= & f(\gamma(t))\,-\, \sum_{k=0}^{n-1} 2^k
\Big[ f(\gamma( 2^{-k-1}  t)^{-2}\gamma( 2^{-k} t))\hspace*{3cm} \label{ncl}\\
& & \qquad\qquad\qquad \quad  +\;
R\big(f(\gamma( 2^{-k-1} t),\,
f(\gamma( 2^{-k-1}t)^{-2}\gamma(2^{-k}t)\big)\Big]\nonumber
\end{eqnarray}
for all $n\in \N$.
Since $\frac{d}{ds}\big|_{s=0} \big(\gamma(\half s)^{-2}\gamma(s)\big)
=0$, we have $\gamma(2^{-k-1}t)^{-2}\gamma(2^{-k}t)=
\cO(2^{-2k})$ here and hence
$f(\gamma(2^{-k-1}t)^{-2}\gamma(2^{-k}t))=\cO(2^{-2\alpha k})$,
as $f$ is $h_\alpha$.
Thus $2^k f(\gamma(2^{-k-1}t)^{-2}\gamma(2^{-k}t))=\cO(2^{-(2\alpha-1) k})$.
Likewise, using that first order Taylor remainders
are at most quadratic in the size of
their argument (Lemma~\ref{taylorR}),
we see that
$2^k R\big(f(\gamma( 2^{-k-1} t),\,
f(\gamma( 2^{-k-1}t)^{-2}\gamma(2^{-k}t)\big)=\cO(2^{-(2\alpha-1) k})$
as well. If $\alpha>\half$,
the preceding estimates entail that
the partial sums of the
series in (\ref{ncl})
form a Mackey-Cauchy sequence,
which converges in $L(H)$
as the latter is assumed Mackey complete.
Therefore $\lim_{n\to\infty}
\frac{f(\gamma(2^{-n}t))}{2^{-n}}$
exists in $L(H)$,
and hence so does
$\lambda:=\lim_{n\to\infty} \frac{f(\gamma(2^{-n}t))}{2^{-n}t}$.
Now clearly $\lambda$ gives us a candidate
for $(f\circ \gamma)'(0)$.
Of course, this rough outline
has to be made more precise.
Furthermore, a lot of work remains:
\begin{itemize}
\item
One has to check that $\lambda$ is independent
of~$t$, and that $(f\circ \gamma)'(0)$ really
exists.
\item
It has to be shown
that $(f\circ \gamma)'(0)$
only depends
on $\gamma'(0)$, and that the mapping
$L(G)\to L(H)$,
$\gamma'(0)\mto (f\circ \gamma)'(0)$
is bounded linear.
\item
The general case $\alpha\in \;]0,1]$
has to be reduced to the case where $\alpha> \half$.
\end{itemize}
To master
the last and penultimate
task, it is essential to
discuss not only a single curve~$\gamma$,
but a whole family of curves
$\eta_s:=\eta(s,\sbull)$ for a smooth map $\eta\!: \R^2\to U$.
Therefore most of the actual proof in Section~\ref{hardlabor}
is formulated for $\eta$'s
instead of mere $\gamma$'s.\\[3mm]
{\bf Organization of the paper.}
In Section~\ref{secprel}, we recall several
basic definitions from convenient differential calculus,
explain some notations,
and compile and develop
various basic facts for later use.
We then define conveniently
H\"{o}lder ($h_\alpha$-) maps
and characterize
them by their behaviour
on ``bornologically compact''
sets (Section~\ref{sechoeld}),
along the lines of the Lipschitz case
treated in~\cite{FaK}
(cf.\ also \cite{FaF}).
In Sections~\ref{secdpoin} and~\ref{secpointw},
we specify and discuss notions of
differentiability at a point
for curves and general mappings, respectively.
H\"{o}lder differentiable curves
and the corresponding mappings
are defined and discussed in
Section~\ref{sechdiff},
as far as required for our purposes.\footnote{See
\cite{FaF} for the general theory
of such maps, which
parallels the familiar $\Lip^k$-case
as in \cite{FaK} or \cite{KaM}.}
In Section~\ref{scpw}, we carry out Step~1
of the proof of our Main Theorem
(curve differentiability at~$1$
implies smoothness).
Before we can carry out Step~2,
further preparations are necessary:
To enable the reduction from arbitrary
$\alpha$ to $\alpha>\half$,
we characterize $h_\alpha$-homomorphisms
(Section~\ref{sectestha})
and in Section~\ref{secrems},
we study the behavior
of Taylor remainders on bornologically compact sets.
In Section~\ref{hardlabor},
the core of the article,
we then complete the proof of our Main Theorem,
based on the reduction steps and preparatory considerations
carried out before.
Various proofs
(part of which are mere adaptations
of the Lipschitz case, \cite{FaK})
have been relegated to
an appendix, and can be taken on faith on a first reading.
\section{Preliminaries}\label{secprel}
This article
is based on
the Convenient Differential
Calculus of Fr\"{o}licher, Kriegl and Michor,
and we presume
familiarity with its basic ideas.
Our main references are \cite{FaK} and
\cite{KaM}.
For the readers convenience, we briefly
recall some of the basic concepts now,
and explain
our notation and terminology.
We also prove various simple results,
for later use.
\begin{numba}\label{defnEB}
Given a locally convex (Hausdorff real topological
vector) space~$E$
and absolutely convex, bounded
subset $B\not=\emptyset$ of $E$,
we let $E_B:=\Spann(B)\sub E$
and make $E_B$ a normed space
with the Minkowski functional $\|.\|_B\!: E_B\to [0,\infty[$,
$\|x\|_B:=\inf\{r>0\!: x\in r B\}$
as the norm.
Then the inclusion map $j\!: E_B\to E$
is continuous linear
(see \cite[Ch.\,III, \S1, No.\,5]{Bou}
for further information).
A locally convex space $E$ is called
a {\em convenient vector space\/} (or: {\em Mackey complete\/})
if $E_B$ is complete,
for each absolutely convex, closed, bounded subset $B\not=\emptyset$
of~$E$ (see \cite[Thm.\,2.14]{KaM}
for alternative characterizations).
\end{numba}
\begin{numba}
The {\em $c^\infty$-topology\/}
on a locally convex space~$E$ is the final topology on~$E$
with respect to the set of all smooth curves
$\gamma\!: \R\to E$ (which are defined as expected).
We write $c^\infty(E)$ for $E$, equipped with the
$c^\infty$-topology.
A subset $U\sub E$ is called {\em $c^\infty$-open\/}
if it is open in $c^\infty(E)$;
in this case, we write $c^\infty(U)$
for $U$, equipped with the topology
induced by~$c^\infty(E)$.
We recall that the $c^\infty$-topology
is finer than the locally convex topology
and can be properly finer;
if $E$ is metrizable, then $c^\infty(E)=E$.
The $c^\infty$-topology
on a product $E\times F$ is finer then
the product topology on $c^\infty(E)\times c^\infty(F)$,
and can be properly finer.
\end{numba}
\begin{numba}
Let $E$ and $F$ be convenient vector
spaces
and $f\!: U\to F$ be a map,
where $U\sub E$ is
$c^\infty$-open. We call $f$ {\em conveniently smooth\/}
or a {\em
$c^\infty$-map\/} if $f\circ \gamma\!: \R\to F$
is a smooth curve, for each smooth
curve $\gamma\!: \R\to E$ with image in~$U$.
If $f\!: U\to F$ is $c^\infty$,
then
the iterated
directional derivatives
$d^kf(x,y_1,\ldots,y_k):=(D_{y_1}\cdots D_{y_k}f)(x)$
exist for any $(x,y_1,\ldots, y_k)\in U\times E^k$
and define a
$c^\infty$-map $d^kf\!: U\times E^k\to F$.
For each $x\in U$, $f^{(k)}(x):=d^kf(x,\sbull)\!: E^k\to F$
is a bounded, symmetric, $k$-linear mapping, and
the map $f^{(k)}\!: U\to \cL^k(E,F)$
into the space of such mappings
(equipped with its natural convenient vector topology)
is $c^\infty$. We abbreviate $df:=d^1f$, $f':=f^{(1)}$,
$f'':=f^{(2)}$.
\end{numba}
\begin{la}\label{EBincl}
Let $E$ be a convenient vector space
and $B\sub E$ be a non-empty,
absolutely convex, bounded subset.
Then the inclusion map $j\!: E_B\to E$
is continuous as a map into~$E$,
equipped with the $c^\infty$-topology.
\end{la}
\begin{proof}
As recalled in {\bf \ref{defnEB}},
$j$ is a continuous linear map into $E$,
equipped with its locally convex vector topology.
Hence $j$ is a bounded linear map
and hence $c^\infty$ \cite[Cor.\,2.11]{KaM}.
Being a $c^\infty$-map, $j$ is continuous
as a map from $c^\infty(E_B)$ to $c^\infty(E)$
(this is immediate from the definition of the
$c^\infty$-topologies).
Here $c^\infty(E_B)=E_B$ due the metrizability of~$E_B$\linebreak
\cite[Thm.\,4.11\,(1)]{KaM}. The assertion follows.
\end{proof}
\begin{numba}
The manifolds and Lie groups of convenient differential calculus
based on the above $c^\infty$-maps
(modelled on convenient vector spaces)
will be referred to as {\em $c^\infty$-manifolds},
resp., {\em $c^\infty$-Lie groups\/} in this article.
The {\em $c^\infty$-topology\/} (or ``natural topology'')
on a $c^\infty$-manifold~$M$ is defined
as the final topology with respect to the set
of smooth curves in~$M$ (see \cite[\S27.4]{KaM}
for further information).
If $f\!: M\to E$ is a $c^\infty$-map from
a $c^\infty$-manifold to a convenient vector space~$E$,
identifying $TE$ with $E\times E$
the tangent map attains the form
$Tf=(f,df)$ for a unique $c^\infty$-map $TM\to E$ denoted~$df$.
\end{numba}
The next lemma is a variant of \cite[Cor.\,2.11]{KaM}:
\begin{la}\label{bdlin}
Let $\alpha\!: E\to F$ be a linear map between
locally convex spaces.
If $\alpha\circ \gamma\!:\R\to F$
is continuous at $0$ for each smooth curve $\gamma\!: \R\to E$
such that $\gamma(0)=0$,
then $\alpha$ is bounded.
\end{la}
\begin{proof}
The proof is by contraposition.
If the linear map~$\alpha$ is not bounded, then there exists
a bounded subset $X\sub E$ such that $\alpha(X)\sub F$
is not bounded. Thus $\|\alpha(X)\|_q$
is unbounded for some continuous seminorm
$\|.\|_q$ on~$F$, entailing that there exist
elements $x_n\in X$ such that $\|\alpha(x_n)\|_q\geq n2^n$.
Then the sequence $(2^{-n}x_n)_{n\in\N}$ in~$E$
converges fast to~$0$ (in the sense of \cite[\S2.8]{KaM}).
Hence, by the Special Curve Lemma \cite[\S2.8]{KaM},
there is a smooth curve
$\gamma\!: \R\to E$ such that $\gamma(\frac{1}{n})=2^{-n}x_n$
for each $n\in \N$.
Since $\|\alpha(\gamma(\frac{1}{n}))\|_q
=2^{-n}\|\alpha(x_n)\|_q\geq n$ tends to $\infty$ as $n\to 0$,
the curve $\alpha\circ \gamma$ is discontinuous at~$0$.
\end{proof}
The following useful fact is clear from the (somewhat
sketchy) discussions
in \cite[end of \S7]{BGN}. For the convenience of the reader,
a self-contained proof is given~in Appendix~A.
\begin{la}\label{BGNcalc}
Let $E$ and $F$ be convenient vector spaces
and $f\!: U\to F$ be a $c^\infty$-map
on a $c^\infty$-open subset $U\sub E$.
Then $U^{[1]}:=\{(x,y,t)\in U\times E\times\R\!:
x+ty\in U\}$ is $c^\infty$-open in $E\times E\times \R$,
and the following map is $c^\infty$:
\[
f^{[1]}\!: U^{[1]}\to F,\qquad
f^{[1]}(x,y,t):=\left\{
\begin{array}{cl}
\frac{f(x+ty)-f(x)}{t} & \mbox{if $\;t\not=0\,${\rm ;}}\\
df(x,y) & \mbox{if $\;t=0\,.$}
\end{array}
\right.
\]
\\[-10mm]\Punkt
\end{la}
\section{Conveniently H\"{o}lder maps}\label{sechoeld}
In this section, we define
and study conveniently H\"{o}lder
mappings. These generalize
the $\cL ip^0$-maps familiar
from convenient differential calculus
(see~\cite{FaK} or \cite{KaM}),
and can be discussed along similar lines
(see also \cite{FaF}).
Throughout the following,
$\alpha\in \;]0,1]$.\\[3mm]
We begin with the definition
of H\"{o}lder
continuous maps on subsets of normed spaces,
and describe their basic
properties, for later use.
\begin{defn}\label{hoelnormed}
Let $(E,\|.\|)$ be a normed real vector space,
$F$ be a real locally convex space, and
$U\sub E$.
A map $f\!: U\to F$ is called {\em H\"{o}lder continuous
of exponent~$\alpha$\/}
(or $H_\alpha$, for short)
if, for every $x\in U$
and continuous seminorm $\|.\|_q\!: F\to [0,\infty[$ on~$F$,
there exists a neighbourhood $V$ of~$x$ in~$U$
and $C\in [0,\infty[$ such that
\begin{equation}
\|f(z)-f(y)\|_q\leq C \, \|z-y\|^\alpha\quad
\mbox{for all $y,z\in V$.}
\end{equation}
$H_1$-maps are also called {\em Lipschitz continuous.}
\end{defn}
Compare \cite[App.\,B]{HOL} for
the case of open domains in locally convex spaces.
\begin{la}\label{laobvio}
Let $\alpha,\beta\in\;]0,1]$.
\begin{itemize}
\item[\rm (a)]
If $\alpha \geq \beta$,
then any $H_\alpha$-map is also $H_\beta$.
\item[\rm (b)]
If $E_1,E_2$ are normed spaces, $U_1\sub E_1$, $U_2\sub E_2$,
$f\!: U_1\to U_2$ is $H_\alpha$ and
$g\!: U_2\to F$ an $H_\beta$-map into a locally
convex space~$F$, then $g\circ f\!:
U_1\to F$ is $H_{\alpha\beta}$.
\item[\rm (c)]
Every smooth curve $\gamma\!: \R\to E$ in a locally
convex space~$E$ is $H_\alpha$.
\end{itemize}
\end{la}
\begin{proof}
(a) and (b) are obvious.
Every smooth curve being Lipschitz continuous
(\cite[\S1.2]{KaM} or \cite[La.\,B2]{HOL}),
(c) follows from (a).
\end{proof}
The special behaviour of H\"{o}lder
maps on compact sets will be exploited
extensively.
\begin{la}\label{hoeloncp}
Let $(E,\|.\|)$ be a
normed space, $K\sub E$ be compact
and $f\!: K\to F$ be an $H_\alpha$-map into
a real locally convex space~$F$.
Then the following holds:
\begin{itemize}
\item[\rm (a)]
$W_K:=\left\{
\frac{f(z)-f(y)}{\|z-y\|^\alpha}\!: y,z\in K,\, y\not=z\right\}$
is bounded in~$F$.
\item[\rm (b)]
There exists an absolutely convex,
bounded subset $D \not=\emptyset$
of $F$
such that $f(K)\sub F_D$,
$W_K\sub D$
and $\|f(z)-f(y)\|_D\leq \|z-y\|^\alpha$
for all $y,z\in K$. In particular,
$W_K$ is bounded
in~$F_D$, and $f \! : K\to F_D$
is~$H_\alpha$ $($and thus continuous$)$.
\end{itemize}
\end{la}
\begin{proof}
(a) If we can show that $\lambda(W_K)\sub \R$
is bounded for each continuous linear functional
$\lambda\!: F\to\R$,
then $W_K$ is bounded by Mackey's Theorem;
here $\lambda\circ f$ is $H_\alpha$
(as is readily verified).
Hence $F=\R$ without loss of generality.
Since $f$ is $H_\alpha$ and
$K$ is compact, we find $n\in \N$,
$C>0$ and
open subsets $V_1,\ldots, V_n,\wt{V}_1,\ldots, \wt{V}_n$ of~$K$
such that $K=\bigcup_{j=1}^n V_j$,
$\wb{V_j}\sub \wt{V}_j$
for $j\in \{1,\ldots, n\}$
(where $\wb{V_j}$ is
the closure of $V_j$ in~$K$),
and
\[
|f(z)-f(y)|\leq C\, \|z-y\|^\alpha\quad\mbox{for all $j\in \{1,\ldots, n\}$
and $y,z\in \wt{V}_j$.}
\]
Thus $X_j:=\{\|z-y\|^{-\alpha}(f(z)-f(y))\!: y,z\in \wt{V}_j,\, y\not=z\}$
is bounded for each $j\in \{1,\ldots, n\}$.
Also $Y_j:=\{\|z-y\|^{-\alpha}(f(z)-f(y))\!: y\in \wb{V_j},\, z\in
K\setminus \wt{V}_j\}$ is bounded,
the sets $K_j:=K\setminus \wt{V}_j$ and
$\wb{V_j}$ being compact and disjoint.
Hence
$W_K\sub \bigcup_{j=1}^n X_j\, \cup\,
\bigcup_{j=1}^n (Y_j\cup (-Y_j))$
is bounded.\vspace{2mm}

(b) We may assume that $K\not=\emptyset$.
The set $W_K$ being bounded by (a),
also the absolutely convex hull~$D$
of $W_K\cup f(K)$ is bounded, and clearly $D$
has the desired properties.
\end{proof}
\begin{la}\label{hoelvsbd}
Let $(E,\|.\|)$ be a finite-dimensional
normed vector space, $U\sub E$ be open,
and $f\!: U\to F$ be a map into
a real locally convex space~$F$.
Then the following holds:
\begin{itemize}
\item[\rm (a)]
$f$ is $H_\alpha$ if and only if
every $x\in U$ has a neighbourhood $K$
for which $W_K\sub F$ $($as in Lemma~{\rm \ref{hoeloncp}\/}$)$
is bounded.
\item[\rm (b)]
$f$ is $H_\alpha$ if and only if $\lambda\circ f$
is $H_\alpha$ for each continuous linear functional
$\lambda\!: E\to \R$. 
\item[\rm (c)]
If $f$ is $H_\alpha$
and $K\sub U$ is compact, then
there exists an
absolutely convex, bounded subset $B\not=\emptyset$
of~$F$ such that $f(K)\sub F_B$
and $f|_K\!: K\to F_B$ is $H_\alpha$.
\item[\rm (d)]
If $f$ is $H_\alpha$, then $f$ is continuous
as a map into $c^\infty(F)$.
\end{itemize}
\end{la}
\begin{proof}
(a) If $f$ is $H_\alpha$ and $x\in U$,
then $W_K$ is bounded for every compact neighbourhood
$K\sub U$ of $x$, by Lemma~\ref{hoeloncp}.
Conversely, assume that every $x\in U$ has a neighbourhood
$K\sub U$ such that $W_K$ is bounded.
Let $\|.\|_q$ be a continuous seminorm
on~$F$. Then
$\|W_K\|_q\sub [0,C]$ for some $C\in \;[0,\infty[$
and thus
$\|f(z)-f(y)\|_q\leq C\, \|z-y\|$ for all $y,z\in K$,
showing that $f$ is $H_\alpha$.\vspace{1mm}

(b) If any $\lambda\circ f$ is $H_\alpha$,
then $W_K$ is weakly bounded and hence bounded,
for any compact neighbourhood
$K\sub U$ of a given element $x\in U$.
Thus $f$ is $H_\alpha$.
The converse is trivial.\vspace{1mm}

(c) is a special case of Lemma~\ref{hoeloncp}\,(b).\vspace{1mm}

(d) Given $x\in U$, we choose a compact
neighbourhood $K\sub U$ of~$x$.
Then $f|_K^{F_B}$ is continuous for $B$ as in\,(c),
and hence so is $f|_K^{c^\infty(F)}$,
by Lemma~\ref{EBincl}.
\end{proof}
\begin{rem}\label{bil1}
Lemma~\ref{hoelvsbd}\,(a)
becomes false for every infinite-dimensional
normed space $(E,\|.\|)$.
Indeed, for any such~$E$, there exists a smooth map
$g\!: E\to\R$ which is unbounded
on the unit ball $B_1(0)\sub E$ (see \cite[La.\,2.3]{Bil}).
Then also
$f\!: E\to \R^\N$, $f(x):=(g(nx))_{n\in \N}$
is smooth
and hence Lipschitz continuous
in the sense of Definition~\ref{hoelnormed},
by \cite[La.\,B2\,(e)]{HOL}.
By construction, $f$
is unbounded on any $0$-neighbourhood~$K\sub E$.
Hence also
$W_K:=\{\|z-y\|^{-1}(f(z)-f(y))\!: 
y,z\in K, y\not=z\}\sub \R^\N$ is unbounded for any~$K$.
\end{rem}
\begin{defn}\label{defconhoel}
Let $E$ be a convenient vector space.
A map $f\!: U\to F$ from a $c^\infty$-open subset $U\sub E$
to a convenient vector space~$F$
is called {\em conveniently H\"{o}lder with exponent~$\alpha$\/}
(or an {\em $h_\alpha$-map\/}, for short)
if $f\circ \gamma$ is an $H_\alpha$-curve,
for each smooth curve $\gamma\!: \R\to U$.
The $h_1$-maps will also
be called {\em conveniently Lipschitz\/}
(or $\Lip^0$, for short).\footnote{In \cite{KaM},
the notation $\cL ip^0$ is used.}
\end{defn}
\begin{rem}\label{whoeldhoeld}
For $f\!: E\supseteq U\to F$
as before, we have:
\begin{itemize}
\item[(a)]
$f$
is $h_\alpha$ if and only if $\lambda\circ
f$ is $h_\alpha$ for each continuous linear functional
$\lambda$ on~$F$, as a consequence
of Lemma~\ref{hoelvsbd}\,(b).
\item[(b)]
If $f$ is $h_\alpha$ and $\gamma\!: I\to U$
is a $c^\infty$-curve defined on an open subset
$I\sub \R$, then $f\circ \gamma$ is~$H_\alpha$.
Indeed, given $t_0\in I$,
a smooth cut-off function $\chi\!:\R\to I$
can be used to create a smooth curve
$\gamma\circ \chi\!: \R\to U$
which coincides with $\gamma$ on a neighbourhood~$J$
of~$t_0$. Then $f\circ \gamma|_J=f\circ (\gamma\circ \chi)|_J$
is $H_\alpha$.
Trivial facts like this one
will be used frequently
in the following, without mention.
\end{itemize}
\end{rem}
We record some immediate
consequences of Lemma~\ref{laobvio}:
\begin{la}\label{trivi}
In the situation of Definition {\rm \ref{defconhoel}},
we have:
\begin{itemize}
\item[\rm (a)]
If $f$ is $h_\alpha$, then $f$
is $h_\beta$ for each $\beta\in \;]0,\alpha]$.
\item[\rm (b)]
If $f$ is $c^\infty$, then $f$ is
$h_\alpha$
for each $\alpha\in \;]0,1]$.\Punkt
\end{itemize}
\end{la}
The following lemma
is a variant
of \cite[Thm.\,1.4.2 and Thm.\,4.3.8]{FaK},
and can be proved analogously.
For completeness, the proof
is given in Appendix~\ref{appB}.
Recall that a subset $K$ of a locally convex space~$E$
is called {\em bornologically compact\/}
if there exists an absolutely convex, bounded
subset $B\not=\emptyset$ of~$E$
such that $K\sub E_B$ and $K$ is compact in~$E_B$.
It is essential
for our purposes that
$h_\alpha$-maps are $H_\alpha$
on bornologically compact sets
(and are, indeed, characterized by this property):
\begin{la}\label{onbcompact}
Let $E$ and $F$ be convenient vector spaces,
$U\sub E$ be $c^\infty$-open, and
$f\!: U\to F$ be a map.
Then the following conditions are equivalent:
\begin{itemize}
\item[\rm (a)]
$f$ is $h_\alpha$.
\item[\rm (b)]
For each absolutely convex, bounded subset
$B\not=\emptyset$ of~$E$ and
compact set $K\sub E_B\cap U$,
the map
$f|_K\!: E_B\supseteq
K\to F$ is $H_\alpha$
$($and thus Lemma~{\rm \ref{hoeloncp}} applies to $f|_K)$.\Punkt
\end{itemize}
\end{la}
For $f\!: \R^n\to\R$,
the next lemma is
covered by \cite[Thm.\,2]{Bom}.
\begin{la}\label{repairKM}
Let $E$ be a finite-dimensional
normed vector space, $U\sub E$ be open,
and $f\!: U\to F$
be a map into a convenient vector space~$F$.
Then $f$ is $h_\alpha$ if and only if $f$ is $H_\alpha$.
\end{la}
\begin{proof}
If $f$ is $H_\alpha$,
then $f\circ \gamma$ is $H_\alpha$
for each smooth curve $\gamma\!: \R\to U$
(Lemma~\ref{laobvio}),
and thus $f$ is $h_\alpha$.
If, conversely,
$f$ is $h_\alpha$,
using that
every $x\in U$ has a compact neighbourhood,
we deduce from Lemmas~\ref{onbcompact},
\ref{hoeloncp} and
\ref{hoelvsbd}\,(a)
that $f$ is $H_\alpha$.
\end{proof}
\begin{la}\label{comphoelder}
Let $E$, $F$ and $H$ be convenient vector
spaces, $U\sub E$ and $V\sub F$ be $c^\infty$-open,
$f\!: U\to V$ be an $h_\alpha$-map,
and $g\!: V\to H$ be an $h_\beta$-map, where
$\alpha,\beta\in \;]0,1]$.
Then
$g\circ f\!: U\to H$ is~$h_{\alpha\beta}$.
\end{la}
\begin{proof}
Let $\gamma\!: \R\to U$ be a smooth curve.
Then $\eta:=f\circ \gamma\!: \R\to V$ is $H_\alpha$.
Given $t_0\in \R$,
set $I:=[t_0-1,t_0+1]$;
there exists an absolutely convex, bounded
subset $B\not=\emptyset$ of~$F$
such that $\eta(I)\sub F_B$
and $\eta|_I\!: I\to F_B$
is $H_\alpha$ (Lemma~\ref{hoelvsbd}\,(c)).
Hence $K:=\eta(I)\sub V$
is compact in $F_B$.
By Lemma~\ref{hoeloncp}\,(b),
there exists an absolutely convex,
bounded subset $D\not=\emptyset$
of~$H$ such that
$g(K)\sub H_D$ and $g|_K\!: F_B\supseteq K\to H_D$
is $H_\beta$. By Lemma~\ref{laobvio}\,(b),
the composition $g|_K\circ (f\circ \eta|_I)$
is $H_{\alpha\beta}$ as a map into $H_D$.
The inclusion map $H_D\to H$ being continuous linear
and thus $H_1$, we deduce that $g\circ (f\circ \eta|_I)$
is $H_{\alpha\beta}$ also as a map into $H$,
using Lemma~\ref{laobvio}\,(b).
Thus $g\circ f\circ \gamma$ is locally $H_{\alpha\beta}$
and hence $H_{\alpha\beta}$. Hence $g\circ f$ is $h_{\alpha\beta}$.
\end{proof}
\begin{la}\label{welldef}
Let $E_j$ be convenient vector spaces
for $j\in \{1,2,3,4\}$,
and $U_j\sub E_j$ be $c^\infty$-open subsets
for $j\in \{1,2,3\}$.
Let
$\phi\!: U_1\to U_2$
and $\psi\!: U_3\to E_4$
be $c^\infty$-maps, and
$f\!: U_2\to U_3$ be $h_\alpha$.
Then $\psi\circ f\circ \phi\!: U_1\to E_4$
is $h_\alpha$.
\end{la}
\begin{proof}
$c^\infty$-maps being $h_1$ by Lemma~\ref{trivi}\,(b),
the assertion follows from Lemma~\ref{comphoelder}.
\end{proof}
\begin{defn}
Let $M$ be a $c^\infty$-manifold modelled
on a convenient vector space~$E$.
A curve $\gamma\!: I\to M$
is called $H_\alpha$
if it is continuous
with respect to
the $c^\infty$-topology on~$M$
and,
for every $t\in I$, there exists a chart
$\phi\!: U\to V\sub E$ of $M$
around $\gamma(t)$
such that $\phi\circ \gamma|_{\gamma^{-1}(U)}\!: \gamma^{-1}(U)\to E$
is $H_\alpha$.
This then holds for all charts,
by Lemma~\ref{welldef}.
\end{defn}
\begin{defn}
A map $f\!: M\to N$ between $c^\infty$-manifolds
is {\em conveniently H\"{o}lder
with exponent $\alpha$\/}
(or an {\em $h_\alpha$-map\/})
if $f\circ \gamma\!: \R\to N$ is an $H_\alpha$-curve
for each smooth curve $\gamma\!: \R\to M$.
If $f$ is $h_\alpha$ for some $\alpha$,
we say that $f$ is {\em conveniently H\"{o}lder}.
\end{defn}
\begin{rem}\label{hoeldthencts}
Any $h_\alpha$-map $f\!: M\to N$ is continuous
with respect to the $c^\infty$-topologies
on~$M$ and~$N$.
Indeed, the topology on~$M$ being final with
respect to the set of smooth curves in~$M$,
this follows from the observation that
$f\circ \gamma$ is an
$H_\alpha$-curve and hence continuous with respect
to the $c^\infty$-topology on~$N$,
for each smooth curve $\gamma\!: \R\to M$.
\end{rem}
\section{Differentiability properties of curves at a point}\label{secdpoin}
In this section, we fix our terminology
concerning differentiability
of curves at a given point
and prove a variant of the
Chain Rule, which enables us to give a meaning
to pointwise differentiability of a curve with
values in a manifold.
\begin{defn}\label{difpropc}
Let $E$ be a convenient vector space,
$\gamma\!: I\to E$ be a map, defined on a
subset
$I\sub \R$, and $t\in I$ be a cluster point of~$I$.
\begin{itemize}
\item[(a)]
$\gamma$ is called {\em differentiable at $t$\/} if
\begin{equation}\label{trivia}
\gamma'(t):=\lim_{s\to t} \frac{\gamma(s)-\gamma(t)}{s-t}
\end{equation}
(where $s\in I\setminus\{t\}$)
exists in~$E$, equipped with its locally convex
vector topology.
\item[(b)]
We say that $\gamma$ is {\em bornologically
differentiable at $t$\/} if there exists a neighbourhood
$J\sub I$ of~$t$ and an absolutely convex, bounded
subset $B\not=\emptyset$ of~$E$
such that $\gamma(J)\sub E_B$ and the limit (\ref{trivia})
exists in $E_B$.\,\footnote{In other words,
the net of difference quotients is Mackey
convergent to $\gamma'(t)$, cf.\ \cite[La.\,1.6]{KaM}.}
\end{itemize}
\end{defn}
If $\gamma$ is bornologically differentiable
at~$t$,
then apparently $\gamma$ is continuous at~$t$
as a map into $E_B$,
and $\gamma$ is differentiable at~$t$
(cf.\ {\bf \ref{defnEB}}).
Beside the standard case
where $I$ is an interval,
we shall encounter the situation where
$t=0$ and $I=\{0\}\cup \{2^{-n}\!: n\in \N\}$.
\begin{la}\label{bornolgood}
Suppose that $\gamma\!: I\to E$ as before is
bornologically differentiable at~$t$ and $f\!: U\to F$
a $c^\infty$-map from an open neighbourhood
$U\sub c^\infty(E)$ of $\gamma(t)$
to a convenient vector space~$F$.
Then $J:=\gamma^{-1}(U)$ is a neighbourhood
of $t$ in~$I$, and $f\circ \gamma|_J\!: J\to F$
is differentiable at~$t$, with derivative
$(f\circ \gamma)'(t)=f'(\gamma(t)).\gamma'(t)$.
\end{la}
\begin{proof}
Choose a neighbourhood
$J\sub I$ of~$t$ and an
absolutely convex, bounded subset $B\not=\emptyset$
of~$E$ as in Definition~\ref{difpropc}\,(b).
Thus $\gamma(J)\sub E_B$, and $\xi\!: J\to E_B$,
\[
\xi(s)\, :=\, 
\left\{
\begin{array}{cl}
\frac{\gamma(s)-\gamma(t)}{s-t} & \mbox{if $\,s\in J\setminus\{t\}$}\\
\gamma'(t) & \mbox{if $\;s=t$}
\end{array}
\right.
\]
is continuous at~$t$ as a map into~$E_B$.
This entails that
$\gamma|_J$ is continuous at~$t$ as a map into
$E_B$ and hence also
as a map into $c^\infty(E)$ (Lemma~\ref{EBincl}).
Therefore $\gamma(J)\sub U$
without loss of generality, after shrinking~$J$.
We have
\[
\frac{f(\gamma(s))-f(\gamma(t))}{s-t}
=
\frac{f\Big(\gamma(t)+(s-t)\frac{\gamma(s)-\gamma(t)}{s-t}\Big)
-f(\gamma(t))}{s-t}
=\, f^{[1]}\bigl(\gamma(t),\, \xi(s),\, s-t\bigr)
\]
for all $s\in J\setminus\{t\}$,
where $f^{[1]}\!: U^{[1]}\to F$ is as in Lemma~\ref{BGNcalc}.
Here
$J\to E\times E \times \R$,
$s\mto (\gamma(t),\xi(s),s-t)$
is continuous at $t$ as a map into $E_B\times E_B\times \R=
(E\times E\times \R)_{B\times B\times [-1,1]}$
and hence also continuous at~$t$
as a map into $c^\infty(E\times E\times \R)$,
by Lemma~\ref{EBincl}.
Since $f^{[1]}$ is $c^\infty$ (see Lemma~\ref{BGNcalc})
and hence continuous
with respect to the $c^\infty$-topologies, we infer that
the limit
$(f\circ \gamma)'(t)= \lim_{s\to t}
f^{[1]}(\gamma(t),\xi(s),s-t)
= f^{[1]}(\gamma(t),\xi(t),0)
=f'(\gamma(t)) . \gamma'(t)$
exists in $c^\infty(F)$ and hence also in~$F$,
and satisfies the required identity.
\end{proof}
\begin{defn}\label{definv}
Let $t\in \R$
and $\gamma\!: I\to M$
be a map from a neighbourhood $I\sub \R$
of~$t$ to a $c^\infty$-manifold~$M$
modelled on a convenient vector space~$E$.
We say that
$\gamma$ is {\em invariantly
differentiable at~$t$\/}
if the following conditions are satisfied:
\begin{itemize}
\item[(a)]
$\gamma$ is continuous at~$t$ (with respect
to the $c^\infty$-topology on~$M$).
\item[(b)]
$\phi\circ \gamma\!: \gamma^{-1}(U)\to E$
is differentiable at~$t$,
for every chart $\phi\!: U\to V\sub E$ of~$M$
around $\gamma(t)$.
\item[(c)]
For any charts $\phi$ and $\psi$ of~$M$
around $\gamma(t)$, we have
$(\psi\circ \gamma)'(t)= A.(\phi\circ \gamma)'(t)$
with $A:=(\psi\circ \phi^{-1})'(\phi(\gamma(t)))$.
\end{itemize}
Property (c) ensures that there is a uniquely
determined element $\gamma'(t)\in T_{\gamma(t)}M$
such that $(d\phi)(\gamma'(t))=(\phi\circ\gamma)'(t)$,
for every chart $\phi$ of $M$ around~$\gamma(t)$. 
\end{defn}
\begin{rem}\label{bornolinv}
Let $M$, $E$ and
$\gamma\!: I\to M$ be as before.
If $\gamma$ is continuous at~$t$
and $\phi\circ \gamma$
is bornologically
differentiable at~$t$ for some chart $\phi$ of~$M$
around $\gamma(t)$,
then $\gamma$ is invariantly
differentiable at~$t$, as a consequence
of Lemma~\ref{bornolgood}. Hence $\gamma'(t)\in TM$
makes sense.
\end{rem}
\section{H\"{o}lder differentiable curves
and {\boldmath $h_\alpha^1$}-maps}\label{sechdiff}
In this section, we define $k$-times H\"{o}lder
differentiable ($H_\alpha^k$-) curves and
conveniently H\"{o}lder differentiable ($h_\alpha^1$-)
maps.
Beyond the Lipschitz case,
$H_\alpha^k$-curves
will only play a role
for $k=1$ (or $k=0$).
We therefore refrain from discussing
$H_\alpha^k$-curves (and $h_\alpha^k$-maps)
for general~$k$
(cf.\ \cite{FaF} for this),
and focus on the case $k=1$.
As to
higher order differentiability,
the standard results on
$\Lip^k$-curves and maps
from \cite{FaK} and \cite{KaM}
are sufficient for~us.
\begin{defn}
Let $E$ be a convenient vector space
and $I\sub \R$ be open.
A map $\gamma\!: I\to E$ is called an {\em $H_\alpha^1$-curve\/}
if it is differentiable at each $t\in I$
and both $\gamma$ and $\gamma'\!: I\to E$
are $H_\alpha$.
Recursively,
we say that $\gamma$ is an
$H_\alpha^{k+1}$-curve if $\gamma$ is $H_\alpha^1$
and $\gamma'$ is an $H_\alpha^k$-curve.
The $H_1^k$-curves are also called $\Lip^k$-curves
(cf.\ \cite[p.\,9]{KaM}).
\end{defn}
\begin{la}\label{basicsH1}
Let $\gamma\!: I\to E$ be an $H^1_\alpha$-curve.
Then the following holds:
\begin{itemize}
\item[\rm (a)]
$\gamma$ is bornologically differentiable at each $t\in I$.
\item[\rm (b)]
If $U\sub c^\infty(E)$ is an open neighbourhood of~$\gamma(I)$
and $f\!: U\to F$ a $c^\infty$-map
to a convenient vector space~$F$,
then $f\circ \gamma\!: I\to F$ is an $H_\alpha^1$-curve,
and $(f\circ \gamma)'(t)=df(\gamma(t),\gamma'(t))$.
\item[\rm (c)]
If $\eta\!: \R \to I$ is smooth,
then $\gamma\circ \eta$ is an $H^1_\alpha$-curve.
\end{itemize}
\end{la}
\begin{proof}
(a) (cf.\ \cite[\S1.7]{KaM}).
Given $t\in I$, let $J\sub I$ be a compact interval
with $t$ in its interior. By
Lemma~\ref{hoeloncp}\,(b),
there exists
an absolutely convex, bounded subset
$B\not=\emptyset$ of~$E$
such that $\gamma'(J)\sub E_B$, $\gamma(J)\sub E_B$ and
$\|\gamma'(r)-\gamma'(s)\|_B\leq |r-s|^\alpha$
for all $r,s\in J$.
Then
$\big\|\frac{\gamma(t+s)-\gamma(t)}{s}-\gamma'(t)\big\|_B
=\big\|\int_0^1\big(
\gamma'(t+rs)-\gamma'(t)\big) \, dr\big\|_B
\leq |s|^\alpha$
for each $s\in (J-t)\setminus\{0\}$,
as $\|\gamma'(t+rs)-\gamma'(t)\|_B\leq |rs|^\alpha\leq
|s|^\alpha$ for each $r$.
Thus $\gamma'(t)={\displaystyle \lim_{s\to 0}} \,s^{-1}(\gamma(t+s)-\gamma(t))$
in~$E_B$.

(b) Fix $t\in I$.
Since $\gamma$ is bornologically differentiable at~$t$
(by (a)), Lemma~\ref{bornolgood}
shows that $f\circ \gamma$
is differentiable at~$t$, with $(f\circ \gamma)'=
df\circ (\gamma,\gamma')$.
Note that $f\circ \gamma$
and $(f\circ \gamma)'=df\circ (\gamma,\gamma')$
are $h_\alpha$
and hence $H_\alpha$,
by Lemmas~\ref{welldef} and \ref{repairKM}.
Hence $f\circ \gamma$ is $H_\alpha^1$.\vspace{2mm}

(c) By Lemmas~\ref{welldef} and \ref{repairKM},
$\gamma\circ \eta$ is $H_\alpha$.
Both $\gamma$ and $\eta$ being $C^1$,
also the composition $\gamma\circ \eta$ is $C^1$,
with $(\gamma\circ \eta)'(t)=\eta'(t)\cdot \gamma'(\eta(t))$
(cf.\ \cite[Prop.\,1.12]{RES}).
Thus $(\gamma\circ \eta)'=m\circ (\eta',\gamma'\circ \eta)$
is a composition of the continuous bilinear
(and hence smooth) scalar multiplication map
$m \!: \R\times E\to E$ and an $H_\alpha$-curve,
and thus $(\gamma\circ \eta)'$ is $H_\alpha$
(see Lemmas~\ref{welldef} and \ref{repairKM}).
\end{proof}
\begin{defn}
Let $M$ be a $c^\infty$-manifold modelled
on a convenient vector space~$E$,
and $I\sub \R$ be open.
A map $\gamma\!: I\to M$
is called an $\Lip^k$-curve (resp., an $H_\alpha^1$-curve)
if it is continuous
with respect to
the $c^\infty$-topology on~$M$ and,
for every $t\in I$, there exists a chart
$\phi\!: U\to V\sub E$ of $M$
around $\gamma(t)$
such that $\phi\circ \gamma|_{\gamma^{-1}(U)}\!: \gamma^{-1}(U)\to E$
is $\Lip^k$ (resp., $H_\alpha^1$).
This then holds for all charts,
by \cite[Cor.\,12.9]{KaM} (resp., Lemma~\ref{basicsH1}\,(b)).
\end{defn}
\begin{la}\label{charh1}
A map $\gamma\!: I\to M$
from an open set $I\sub \R$
to a $c^\infty$-manifold
is an $H_\alpha^1$-curve if and only
$\gamma$ is invariantly differentiable
at each $t \in I$ and the curve
$\gamma'\!: I\to TM$~is~$H_\alpha$.
\end{la}
\begin{proof}
If $\gamma$ is an $H_\alpha^1$-curve,
then $\gamma$ is continuous
(see Remark~\ref{hoeldthencts}).
Given $t\in I$,
we let $\phi\!: U\to V\sub E$
be a chart of $M$ around $\gamma(t)$.
Then $\phi\circ \gamma|_J$ is an $H^1_\alpha$-curve
on $J:=\gamma^{-1}(U)$ and thus bornologically
differentiable at~$t$ (see Lemma~\ref{basicsH1}\,(a)),
whence $\gamma$ is invariantly
differentiable at~$t$ (Remark~\ref{bornolinv}).
If, conversely,
$\gamma$ is invariantly
differentiable at each~$t$
and $\gamma'$ an $H_\alpha$-curve,
then $\gamma$ is continuous
and for $t$, $\phi$ and $J$ as before
$(\phi\circ \gamma|_J)'=d\phi\circ \gamma'|_J$
is $H_\alpha$ (by Lemmas~\ref{welldef} and \ref{repairKM}),
entailing that $\phi\circ \gamma|_J$
is an $H_\alpha^1$-curve.
\end{proof}
\begin{defn}
A map $f\!: M \to N$ between $c^\infty$-manifolds
is called $h_\alpha^1$ (resp., $\Lip^k$) if
$f\circ \gamma\!: \R\to N$ is an
$H^1_\alpha$-curve (resp., a $\Lip^k$-curve),
for each smooth curve $\gamma\!: \R\to M$.
\end{defn}
\section{Curve differentiability of maps at a given point}\label{secpointw}
In this section,
we define and discuss differentiability
properties of mappings at a given point.
``Bornological
curve differentiability,''
which we introduce here,
is a bornological variant of the
classical notion of Hadamard differentiability.
It is a well-chosen concept
in the sense that, on the one
hand, bornological
curve differentiability at $1$ is a sufficiently
strong property to ensure
smoothness of a homomorphism.
On the other hand, it is a sufficiently weak
differentiability property, in the sense that
we shall manage to establish it
for conveniently H\"{o}lder homomorphisms.
We also
introduce a notion of ``curve differentiability,''
as a technical tool.
This is a weaker and much more illusive concept.\\[3mm]
Let $E$ and $F$ be convenient vector spaces,
$U\sub E$ be $c^\infty$-open,
$x\in U$,
and $f\!: U\to F$.
\begin{defn}\label{defncd}
$f$ is called
{\em curve differentiable
at $x$\/} if the conditions
(a)--(d) are satisfied:
\begin{itemize}
\item[(a)]
$f\circ \gamma$ is differentiable at~$0$, for every smooth curve
$\gamma\!: \R\to U$ such that $\gamma(0)=x$.
\item[(b)]
There exists a (necessarily unique)
bounded linear map $f'(x)\!: E\to F$
such that $(f\circ \gamma)'(0)=f'(x).\gamma'(0)$,
for every $\gamma$ as in (a).
\item[(c)]
$f$ is continuous at~$x$ with respect to
the $c^\infty$-topologies on~$U$ and $F$.
\item[(d)]
For every convenient
vector space $H$ and $c^\infty$-map
$g\!: V\to H$ defined on an
open neighbourhood~$V$
of $f(x)$ in $c^\infty(F)$, and every smooth curve
$\gamma\!: \R\to f^{-1}(V)$,
the curve $g\circ f\circ \gamma$ is differentiable at~$0$, with
$(g\circ f\circ \gamma)'(0)
=g'(f(x)).f'(x).\gamma'(0)$.
\end{itemize}
We call $f$
{\em bornologically curve differentiable at $x$\/}
if $f\circ \gamma$ is bornologically
differentiable at~$0$ for every $\gamma$ as in (a),
and conditions (b) and (c) are satisfied.
\end{defn}
Of course, if $f$ is curve differentiable at~$x$,
then analogues of (a), (b) and (d) hold if
$\gamma\!: I\to U$ is defined on an
open $0$-neighbourhood $I\sub \R$ only
(cf.\ Remark~\ref{whoeldhoeld}\,(b)).
Also note that (c) holds
if $f$ is conveniently H\"{o}lder, by Remark~\ref{hoeldthencts}.
\begin{rem}\label{rem1cd}
Let $f$ as before be curve differentiable
at~$x$, and $g$ be as in (d).
\begin{itemize}
\item[(a)]
If $h\!: W\to f^{-1}(V)$ is a $c^\infty$-map
defined on a $c^\infty$-open subset
of a convenient vector space
and $w\in W$ such that $h(w)=x$,
then
$(g\circ f\circ h\circ \gamma)'(0)=(g\circ f\circ
(h\circ \gamma))'(0)=g'(f(x)).f'(x).(h\circ \gamma)'(0)=
g'(f(x)).f'(x).h'(w).\gamma'(0)$
for any smooth curve $\gamma\!: \R\to W$
such that $\gamma(0)=w$, exploiting
condition\,(d) from
Definition~\ref{defncd}.
\item[(b)]
For $f,g, h$ and $w$ as before,
we readily deduce from
Definition~\ref{defncd}\,(d),
Part\,(a) of the present remark
and the Chain Rule for $c^\infty$-maps
that $g\circ f\circ h\!: W\to H$
is curve differentiable at~$w$,
with $(g\circ f\circ h)'(w)=g'(f(x))\circ f'(x)\circ
h'(w)$.
\end{itemize}
\end{rem}
From Lemma~\ref{bornolgood}, we immediately deduce:
\begin{la}\label{strongstronger}
If $f\!: U\to F$ $($as before$)$
is bornologically curve differentiable at~$x$,
then $f$ is curve differentiable
at~$x$.\Punkt
\end{la}
\begin{defn}\label{defcdifmfd}
Let $M$ and $N$ be $c^\infty$-manifolds
modelled on $E$, resp., $F$,
and $x\in M$.
A map
$f\!:M\to N$
is called {\em curve differentiable at $x$\/}
(resp., {\em bornologically curve differentiable at $x$\/}),
if there
exists a chart $\phi\!: U_1\to U\sub E$ of $M$ around
$x$ and a chart $\psi\!: V_1\to V\sub F$ of $N$ around
$f(x)$ such that $f(U_1)\sub V_1$ and
$\psi\circ f\circ \phi^{-1}$
is curve differentiable (resp.,
bornologically curve differentiable)
at~$\phi(x)$.
\end{defn}
\begin{rem}\label{chaincurve}
Let $f$ (as before)
be curve differentiable at~$x$.
The the following holds:
\begin{itemize}
\item[(a)]
$f$ is continuous at~$x$ with respect to the $c^\infty$-topologies
on~$M$ and $N$.
\item[(b)]
In the situation
of Definition~\ref{defcdifmfd},
$\psi\circ f\circ \phi^{-1}$
is curve differentiable at~$\phi(x)$
for {\em every\/} choice of charts
(use Remark~\ref{rem1cd}\,(b)).
\item[(c)]
If $\gamma\!: I \to M$ is a $c^\infty$-curve
on an open $0$-neighbourhood $I\sub \R$
such that~\mbox{$\gamma(0)=x$,}
then $f \circ \gamma$ is invariantly differentiable
at~$0$ (in view of Part\,(a) of the present
remark and Part\,(d) of the definition
of curve differentiability, applied
to $f$ in local coordinates).
Furthermore, $(f \circ \gamma)'(0)$
only depends on $\gamma'(0)$ (as a consequence
of Definition~\ref{defncd}\,(b)).
Hence $T_xf\!: T_xM\to T_{f(x)}N$,
$(T_xf).\gamma'(0):=(f\circ\gamma)'(0)$
is well defined.
As a consequence of Definition~\ref{defncd}\,(b),
$T_xf$ is bounded linear.
\item[(d)]
Remark~\ref{rem1cd}\,(b) readily entails:
If
$Y,Z$ are $c^\infty$-manifolds,
$g\!: N\to Z$ as well as $h\!: Y\to M$
are $c^\infty$-maps,
and $y\in Y$ is an element such that $h(y)=x$,
then $g\circ f\circ h\!: Y\to Z$ is curve
differentiable at~$y$, with
$T_y(g\circ f\circ h)=T_{f(x)}g\circ T_xf\circ T_yh$.
\end{itemize}
\end{rem}
\section{$\!\!\!$Pointwise differentiable homomorphisms are smooth}\label{scpw}
We now perform Step~1 of the programme outlined
in the Introduction:
curve differentiability
at $1$ implies smoothness for homomorphisms.
A simple observation will be used:
\begin{la}\label{C1smooth}
Let $f\!: G\to H$ be a $\Lip^1$-homomorphism
between
Lie groups in the sense of convenient differential
calculus.
Then $f$ is a $c^\infty$-map.
\end{la}
\begin{proof}
We show that $f$ is of $\Lip^k$
for each $k\in \N$, by induction.
By hypothesis, $f$ is $\Lip^1$;
it therefore gives rise to a tangent map
$Tf\!: TG\to TH$.
With respect to the left trivializations
of $TG$ and $TH$, the map
$Tf$ corresponds to $f\times L(f)\!: G\times L(G)\to H\times L(H)$.
Here $f$ is $\Lip^k$, and the linear map
$L(f)$ is $\Lip^0$ (cf.\
\cite[Thm.\,12.8]{KaM}),
hence bounded (Lemma~\ref{bdlin})
and thus smooth \cite[Cor.\,2.11]{KaM}.
Hence $Tf$ is $\Lip^k$.
Now $f$ being $\Lip^1$ with $Tf$ a $\Lip^k$-map,
the map $f$ is $\Lip^{k+1}$ (cf.\ \cite[Thm.\,12.8]{KaM}).
\end{proof}
\begin{prop}\label{suffs}
Let $f\!: G\to H$
be a conveniently H\"{o}lder
homomorphism
between $c^\infty$-Lie groups.
If $f$ is
curve differentiable at~$1$,
then $f$ is smooth.
\end{prop}
\begin{proof}
Being conveniently H\"{o}lder,
$f$ is $h_\alpha$ for some $\alpha\in \;]0,1]$.
Given $x\in G$ we have $f=\lambda_{f(x)}^H\circ f\circ
\lambda^G_{x^{-1}}$
(with left translation maps as indicated),
because $f$ is a homomorphism.
The map
$f$ being curve differentiable
at~$1$, using Remark~\ref{chaincurve}\,(d),
we deduce from the latter formula that $f$
is curve differentiable at~$x$, with
\begin{equation}\label{newapproach}
T_xf=
T_1\lambda^H_{f(x)}\circ T_1f\circ T_x\lambda_{x^{-1}}^G\,.
\end{equation}
Now let $\gamma\!: \R\to G$ be a
smooth curve.
Given $t\in \R$,
the curve $f\circ \gamma$ is invariantly
differentiable at~$t$
because $f$ is curve differentiable
at $\gamma(t)$ (cf.\ Remark~\ref{chaincurve}\,(c)).
By
(\ref{newapproach}), we have
$(f\circ \gamma)'(t)=(T_{\gamma(t)}f). \gamma'(t)=
(T_1\lambda^H_{f(\gamma(t))} \circ T_1f \circ T_{\gamma(t)}
\lambda_{\gamma(t)^{-1}}^G). \gamma'(t)$.
Using
the bounded linear (and hence $c^\infty$-) map
$A:=T_1f\!: L(G)\to L(H)$,
the Maurer-Cartan form
$\omega_G\!: TG\to L(G)$,\linebreak
$T_xG\ni v\mto (T_x\lambda^G_{x^{-1}})(v)$
(which is $c^\infty$),
and the $c^\infty$-map
$h\!: H\times L(H)\to TH$
defined via
$h(x,v):= (T_1\lambda^H_x)(v)$,
the preceding formula can be rewritten
in the form
\[
(f\circ \gamma)'(t)\;=\;
h\bigl( f(\gamma(t)),A.\omega_G(\gamma'(t))\bigr)
\qquad\mbox{for all $t\in \R\,$.}
\]
The curve
$\zeta:=A\circ \omega_G\circ \gamma'\!:
\R\to L(H)$ is smooth
and thus $H_\alpha$.
Furthermore,
$\eta:=f\circ \gamma\!: \R\to H$
is $H_\alpha$, as $f$ is $h_\alpha$.
Hence $(\eta,\zeta)\!: \R\to H\times L(H)$
is $H_\alpha$ and thus
$(f\circ \gamma)'\,=\, h\circ (\eta,\zeta)$
is $H_\alpha$, by Lemmas~\ref{repairKM}
and \ref{welldef}.
Now
Lemma~\ref{charh1} shows that $f$ is $h_\alpha^1$
and hence $\Lip^0$, using that $h^1_\alpha$-curves
in $L(H)$ are $C^1$ and hence $\Lip^0$.
We may therefore take
$\alpha=1$ in the preceding considerations,
showing that $f$ is $h_1^1=\Lip^1$
and hence smooth (Lemma~\ref{C1smooth}).
\end{proof}
\section{Testing at \boldmath{$\one$} whether a homomorphism is
\boldmath{$h_\alpha$}}\label{sectestha}
In this section,
we explain how the $h_\alpha$-property
of a homomorphism
can be characterized
by a suitable property at the identity element. 
In the Section~\ref{hardlabor},
this characterization
will be used to show that an $h_\beta$-homomorphism,
where $\beta\in\;]0,\half]$,
is also $h_\alpha$ for $\alpha:=\frac{3}{2}\beta$.
\begin{la}\label{hoeldatpoint}
A homomorphism $f\!: G\to H$
between
$c^\infty$-Lie groups
is $h_\alpha$ if and only
if $f$ is continuous with respect
to the $c^\infty$-topologies
and for every
$c^\infty$-map $\theta\!: \R^2\to G$
with $\theta(s,0)=1$ for all $s\in \R$,
there exists a chart $\phi\!: U\to V\sub L(H)$
of~$H$ around~$1$,
an open $0$-neighbourhood $I \sub \R$
such that $f(\theta(I\times I))\sub U$,
and an absolutely convex, bounded subset
$B\not=\emptyset$ of $L(H)$
such that
\begin{itemize}
\item[\rm (a)]
$\phi(f(\theta(I^2)))\sub L(H)_B$;
\item[\rm (b)]
$\xi:=\phi\circ f\circ \theta|_{I^2} \!: I^2\to L(H)_B$
is continuous;
\item[\rm (c)]
There exists $K\geq 0$ such that
$\|\xi(s,t)-\xi(s,0)\|_B\leq K |t|^\alpha$,
for all $s,t\in I$.
\end{itemize}
This then holds for any choice of the chart~$\phi$.
\end{la}
\begin{proof}
The final assertion
is easily established using
Lemma~\ref{onbcompact};
we omit the details.
The necessity of the condition
is apparent
(cf.\ Remark~\ref{hoeldthencts},
Lemma~\ref{repairKM}, Lemma~\ref{hoeloncp}).
\begin{numba}
Conversely, assume now
that $f$ satisfies
the described condition.
We have to show that
$f\circ \gamma$ is $H_\alpha$,
for each $c^\infty$-curve $\gamma\!: \R\to G$.
This will hold if $f\circ \gamma$ is $H_\alpha$
on some open neighbourhood
of each given $t_0\in \R$.
Because $(f\circ \gamma)(t)=f(\gamma(t_0))f(\gamma(t_0)^{-1}\gamma(t))
=(\lambda_{f(\gamma(t_0))}\circ f\circ \zeta)(t-t_0)$,
where $\zeta\!: \R\to G$, $\zeta(t):=\gamma(t_0)^{-1}\gamma(t+t_0)$
is a $c^\infty$-curve,
after replacing $\gamma$ with $\zeta$
it actually suffices to assume that $t_0=0$
and $\gamma(t_0)=1$.
\end{numba}
\begin{numba}
Pick a chart
$\phi\!: U \to V\sub L(H)$ of~$H$
around~$1$ such that
$\phi(1)=0$.
The group
multiplication $m\!: H\times H\to H$ being~$c^\infty$,
$M_1:=\{(x,y)\in U \times U \!: xy\in U\}$
is open in $H\times H$ (equipped with the $c^\infty$-topology).
Then $M:=(\phi\times \phi)(M_1)$
is $c^\infty$-open in $L(H)\times L(H)$, and
$\mu:=\phi\circ m|_{M_1}^U \circ (\phi^{-1}\times \phi^{-1})|_M^{M_1}
\!: M\to V$ is $c^\infty$.
\end{numba}
\begin{numba}
The $c^\infty$-maps
$\theta\!: \R^2\to G$, $\theta(s,t):=\gamma(s)^{-1}\gamma(s+t)$
and
$\bar{\theta}\!: \R^2\to G$, $\bar{\theta}(s,t):=\gamma(t)$
satisfy $\theta(s,0)=\bar{\theta}(s,0)=1$
for all $s\in \R$. Using the hypotheses on~$f$
(and the final assertion of the lemma),
we find an open $0$-neighbourhood $I\sub \R$,
$K_1\geq 0$
and an absolutely convex, bounded subset $B\not=\emptyset$
of $L(H)$ such that $f(\theta(I^2)), f(\bar{\theta}(I^2))\sub U$,
$\phi(f(\theta(I^2))), \phi(f(\bar{\theta}(I^2)))\sub L(H)_B$,
and both of the maps $\xi:=
\phi\circ f\circ \theta|_{I^2}\!: I^2\to L(H)_B$
and $\bar{\xi}:=
\phi\circ f\circ \bar{\theta}|_{I^2}\!: I^2\to L(H)_B$
are continuous and
$\|\xi(s,t)\|_B\leq K_1|t|^\alpha$
for all $s,t\in I$.
Then also the map $\bar{\xi}\times \xi \!: I^2\times I^2
\to c^\infty(L(H)\times L(H))$
is continuous (cf.\ Lemma~\ref{EBincl}).
After shrinking~$I$, we may hence
assume that $\bar{\xi}(I^2)\times \xi(I^2) \sub M$.
\end{numba}
\begin{numba}
Let $A\sub I$ be a compact $0$-neighbourhood
and $J\sub \R$ be an open $0$-neighbourhood
such that $J-J\sub A$.
We claim that
$\phi\circ f\circ \gamma|_J\!: J\to L(H)$
is~$H_\alpha$ (where $\phi(f(\gamma(t)))=\bar{\xi}(0,t)$).
To see this, let $\|.\|_p$ be a continuous seminorm
on~$L(H)$.
The set $C:= (\bar{\xi}\times \xi)(A^2\times A^2)
\sub (L(H)_B)^2=(L(H)^2)_{B^2}$
is compact. Since $\mu$ is $c^\infty$ and thus $\Lip^0$,
using Lemmas~\ref{onbcompact} and \ref{hoeloncp}\,(a)
we find
$K_2\in [0,\infty[$
such that
$\|\mu(u)-\mu(v)\|_p\leq K_2 \|u-v\|_{B^2}$
for all $u,v\in C$.
Given $r,t\in J\sub A$, we have $s:=r-t\in A$.
Using that $f(\gamma(t))=\phi^{-1}(\bar{\xi}(0,t))$,
$\,f(\gamma(t)^{-1}\gamma(t+s))=\phi^{-1}(\xi(t,s))$
and $f(1)=\phi^{-1}(\xi(t,0))$, we obtain
\begin{eqnarray*}
\phi(f(\gamma(r)))-\phi(f(\gamma(t)))&=&
\phi(f(\gamma(t+s)))-\phi(f(\gamma(t)))\\
&=& \phi(f(\gamma(t))f(\gamma(t)^{-1}\gamma(t+s)))-
\phi(f(\gamma(t))f(1))\\
&=&
\mu(\bar{\xi}(0,t),\xi(t,s))-\mu(\bar{\xi}(0,t),\xi(t,0))\,.
\end{eqnarray*}
We now deduce that
$\|\phi(f(\gamma(r)))-\phi(f(\gamma(t)))\|_p=
\|\mu(\bar{\xi}(0,t),\xi(t,s))-\mu(\bar{\xi}(0,t),\xi(t,0))\|_p
\leq K_2 \|\xi(t,s)-\xi(t,0)\|_B
=K_2 \|\xi(t,s)\|_B
\leq
K_1K_2 |s|^\alpha=K_1K_2 |r-t|^\alpha$.
\end{numba}
Thus $\phi\circ f\circ \gamma|_J\!: J\to L(H)$ is $H_\alpha$
and hence so is $f\circ \gamma|_J$,
which completes the proof.
\end{proof}
\section{Estimates on Taylor remainders}\label{secrems}
We study
the behaviour of first order Taylor remainders
on bornologically compact sets.
\begin{la}\label{taylorR}
Let $E$ and $F$ be convenient vector spaces,
$U\sub E$ be $c^\infty$-open
and $f\!: U\to F$
be $c^\infty$; define the ``first
order Taylor remainder'' $\rho \!: U\times U\to F$
of~$f$ via
$\rho (x,y):=f(y)-f(x)-f'(x).(y-x)$.
Let $C \not=\emptyset$
be an absolutely convex, bounded subset of~$E$
and $K\sub E_C\cap U$ be compact.
Then there exists an absolutely convex, bounded
subset $B\not=\emptyset$ of~$F$ such that $f(K)\sub F_B$,
$\rho(K\times K)\sub F_B$, the restrictions
$f|_K\!: E_C\supseteq K\to F_B$ and
$\rho|_{K^2}\!: (E_C)^2\supseteq
K\times K \to F_B$
are Lipschitz continuous, and
\begin{equation}\label{est2ndord}
\|\rho(x,y)\|_B\leq (\|x-y\|_C)^2\quad \mbox{for all $x,y\in K$.}
\end{equation}
\end{la}
\begin{proof}
After replacing $C$ with its closure in~$E$,
we may assume that $C$ is the closed unit
ball in $(E_C,\|.\|_C)$.
The maps $f$ and  $\rho$ being $c^\infty$ and hence $\Lip^0$,
Lemma~\ref{onbcompact}
provides an absolutely convex, bounded subset
$B_0\not=\emptyset$ of~$F$ such that $f(K)\sub F_{B_0}$,
$\rho(K\times K)\sub F_{B_0}$,
and such that both $f|_K\!: E_C\supseteq K\to F_{B_0}$
and $\rho|_{K^2}\!: (E_C)^2\supseteq K^2\to F_{B_0}$
are Lipschitz continuous.
Since $U\cap E_C$ is open in $E_C$ (see Lemma~\ref{EBincl})
and $K$ compact, we find $\ve>0$ such that
$K+\ve C\sub U$. Then the set
$K_\ve:=\{(x,y)\in K^2\!: \|x-y\|_C\leq \ve\}$ is compact in
$(E_C)^2$, and
$[x,y]:=\{x+t(y-x)\!: t\in [0,1]\}\sub
U$, for any $(x,y)\in K_\ve$.
Hence, by Taylor's Formula \cite[Prop.\,4.4.18]{FaK},
we have
\begin{equation}\label{taylint}
\rho(x,y)=
\int_0^1(1-t)f''(x+t(y-x))(y-x,y-x)\, dt\quad
\mbox{for all $(x,y)\in K_\ve$.}
\end{equation}
The set $\wt{K}:=\bigcup_{(x,y)\in K_\ve}[x,y]\sub U$
is compact in $E_C$.
Now $f''\!: U\to \cL^2(E,F)$
being~$c^\infty$, the image $f''(\wt{K})$
is compact in $\cL^2(E,F)$, equipped with
the $c^\infty$-topology.
Hence $f''(\wt{K})$ is also
compact in
the locally
convex space $\cL^2(E,F)$
and thus bounded.
The trilinear evaluation map
\mbox{$\ev\!: \cL^2(E,F)\!\times\! E\!\times \!E \to F$}
being $c^\infty$ (cf.\ \cite[Cor.\,3.13\,(1)]{KaM})
and thus bounded \cite[La.\,5.5]{KaM},
we see that the image
$\ev(f''(\wt{K})\times C \times C)$
is bounded in $F$ and hence
contained in an absolutely convex,
bounded subset $B\not=\emptyset$ of~$F$:
\begin{equation}\label{refadd}
\ev(f''(\wt{K})\times C \times C)\sub B\,.
\end{equation}
After increasing $B$, we may assume
that $B$ is closed and $B_0\sub B$.
As a consequence of (\ref{refadd}),
we have $d^2f(\wt{K}\times (K-K)\times (K-K))\sub F_B$;
since $\wt{K}$ is compact,
after increasing~$B$ further
we may assume that 
and $d^2f$ is Lipschitz continuous
as a map from $\wt{K}\times (K-K)^2\sub (E_C)^3$ to~$F_B$
(Lemma~\ref{onbcompact}).
This entails that the integrands in (\ref{taylint})
are continuous as maps into
the Banach space~$F_B$.
Hence, the integral also
exists in~$F_B$,
and clearly the $F_B$-valued integral
coincides with
the $F$-valued integral $\rho(x,y)$
(equality can be tested with linear functionals
in~$F'$).
Now \cite[Cor.\,2.6\,(4)]{KaM}
(or \cite[La.\,1.7]{RES})
implies that
\begin{equation}\label{specsmall}
\|\rho(x,y)\|_B\leq
\sup_{t\in [0,1]} \|f''(x+t(y-x))(y-x,y-x)\|_B
\leq \|y-x\|_C^2\quad \mbox{for all $(x,y)\in K_\ve$,}
\end{equation}
where (\ref{refadd}) was used to get the second inequality.
To complete the proof,
note that $D:=\{(x,y)\in K^2\!: \|x-y\|_C\geq \ve\}$
is compact in $(E_C)^2$, whence $\rho(D)$ is compact in $F_{B_0}$
and hence bounded in~$F$. Let $B_1:=\absconv(\ve^{-2}\rho(D))$
be the absolutely
convex hull of $\ve^{-2}\rho(D)$,
and replace $B$ with $\absconv(B\cup B_1)$.
Then (\ref{specsmall})
remains valid, and furthermore
$\|\rho(x,y)\|_B=\ve^2 \|\ve^{-2}\rho(x,y)\|_B\leq
\ve^2\leq \|(x,y)\|_C^2$
for all $(x,y)\in D$. Since $K^2=K_\ve\cup D$,
the preceding estimate and (\ref{specsmall})
show that (\ref{est2ndord}) holds.
\end{proof}
\section{Conveniently H\"{o}lder homomorphisms are smooth}\label{hardlabor}
Having completed all necessary preparations,
we are now ready
to prove the main result.
\begin{thm}\label{main}
Let $f\!: G\to H$ be a homomorphism
between Lie groups in the sense of convenient differential
calculus.
If $f$ is conveniently H\"{o}lder,
then $f$ is a $c^\infty$-map.
In particular, $f$ is $c^\infty$ if $f$ is $\Lip^0$.
\end{thm}
\begin{proof}
By hypothesis, $f$ is $h_\alpha$ for some $\alpha\in \;]0,1]$.
We choose a chart $\phi\!: U_1\to U\sub L(G)$
of $G$ around~$1$
and a chart $\psi\!: V_1\to V\sub L(H)$ of~$H$
around~$1$
such that $f(U_1)\sub V_1$,
$\phi(1)=0$, and $\psi(1)=0$.
Then
$g:=\psi\circ f\circ \phi^{-1}\!: U\to V\sub L(H)$
is $h_\alpha$, and $g(0)=0$.
\begin{numba}
The map $\bar{\sigma}\!: V_1\times V_1\to H$, $\bar{\sigma}(x,y):=x^2y$
being $c^\infty$, the preimage $S_1:=\bar{\sigma}^{-1}(V_1)$
is $c^\infty$-open in $V_1\times V_1$.
Then
$S:=(\psi\times \psi)(S_1)$
is a $c^\infty$-open
$(0,0)$-neighbourhood in
$L(H)\times L(H)$, and
$\sigma:=\psi\circ \bar{\sigma}\circ (\psi^{-1}\times
\psi^{-1})|_S \!: S\to V$ is $c^\infty$.
The first order Taylor expansion of $\sigma$
around $(0,0)$ gives
\begin{equation}\label{sigmaR}
\sigma(x,y)=2x+y+R(x,y)\quad \mbox{for all $(x,y)\in S\sub L(H)\times L(H)$,}
\end{equation}
where $R\!: S\to L(H)$, $R(x,y):=\rho_\sigma((0,0), (x,y))$
with $\rho_\sigma\!: S\times S \to L(H)$
the first order
Taylor remainder of~$\sigma$
(as in Lemma~\ref{taylorR}).
\end{numba}
\begin{numba}
The map $\bar{\tau}\!: U_1\times U_1\to G$, $\bar{\tau}(x,y):=x^{-2}y:=
x^{-1}x^{-1}y$
being $c^\infty$, the set
$W_1:=\bar{\tau}^{-1}(U_1)$
is $c^\infty$-open in $U_1\times U_1$.
Then
$W:=(\phi\times \phi)(W_1)$
is a $c^\infty$-open
$(0,0)$-neighbourhood
in $L(G)\times L(G)$,
and
$\tau := \phi\circ \bar{\tau}\circ (\phi^{-1}\times
\phi^{-1})|_W \!: W\to U$ is $c^\infty$.
We have $\tau(0,0)=0$ and
\begin{equation}\label{tauR}
d\tau((0,0),(u,v))=-2u+v\qquad
\mbox{for all $u,v\in L(G)$.}
\end{equation}
\end{numba}
\begin{numba}\label{makestarsh}
Let
$\ell\in \{1,2\}$,
$\delta\in \;]0,\infty]$
and $\eta \!: B(\delta)\to U\sub L(G)$
be a $c^\infty$-map on
$B(\delta):=\;]{-\delta,\delta}[^\ell\sub \R^\ell$,
such that $\eta(s,0)=0$
for all $s\in \;]{-\delta},\delta[^{\ell-1}\sub \R^{\ell-1}$
(if $\ell=1$, we identify $\R$ with
$\R^0\times \R$ here and in the following,
to unify notation. Thus, the
argument $s$ has to be ignored if $\ell=1$,
and what we require
is $\eta(0)=0$).
Then $\theta:=\phi^{-1}\circ \eta \!: B(\delta)\to U_1\sub G$
is~$c^\infty$.
The map $\zeta\!: B(\delta)\to L(G)\times L(G)$,
$\zeta(s,t):=(\eta(s,\half t),\eta(s,t))$
being $c^\infty$
and hence continuous
into
$c^\infty(L(G)\times L(G))$,
we find $\delta_1\in\;]0,\delta]$
with $\zeta(B(\delta_1))\sub W$.
Then $\kappa\!: B(\delta_1)\to U$,
$\kappa(s,t):=\tau(\eta(s,\half t),\eta(s,t))$
is~$c^\infty$.
\end{numba}
\begin{numba}
The map $\chi\!: B(\delta_1)\to H\times H$, $\chi(s,t):=
\bigl( f(\theta(s,\half t)), \, f(\theta(s,\half t)^{-2}\theta(s,t))\bigr)$
being $h_\alpha$,
we find $\delta_2\in \;]0,\delta_1]$
such that $\chi(B(\delta_2))\sub S_1$.
We define $\omega\!: B(\delta_2)\to S$,
$\omega:=(\psi\times \psi)\circ \chi|_{B(\delta_2)}$.
Then
$\omega(s,t)=
\bigl( g(\eta(s,\half t)),\,
g(\tau(\eta(s,\half t),\eta(s,t)))\bigr)
=\bigl( g(\eta(s,\half t)),\,
g(\kappa(s,t))\bigr)$
and
\begin{equation}\label{facilrewr}
\sigma(\omega(s,t))\,=\,\sigma\bigl(
g(\eta(s,\half t)),\,
g(\tau(\eta(s,\half t),\eta(s,t))) \bigr)\,=\, g(\eta(s,t))\;\;
\mbox{for all $(s,t)\in B(\delta_2)$,}
\end{equation}
as
$\psi^{-1}\bigl( \sigma(g(\eta(s,\half t)),\,
g(\tau(\eta(s,\half t),\eta(s,t)))) \bigr)
= f(\theta(\half s,t))^2 f\big(\theta(s,\half t)^{-2}\theta(s,t)\big)
=f(\theta(s,t))$ $=\psi^{-1}\big( g(\eta(s,t)) \big)$.
\end{numba}
\begin{numba}
Combining (\ref{facilrewr})
and (\ref{sigmaR}) yields
$g(\eta(s,t)) =
2g(\eta(s,\half t))+ g(\kappa(s,t))
+R(\omega(s,t))$, whence
\begin{equation}\label{startit}
g(\eta(s,\half t))\,=\,
\half g(\eta(s,t)) -\half g(\kappa(s,t))
-\half R(\omega(s,t))
\quad \mbox{for all $(s,t)\in B(\delta_2)$.}
\end{equation}
Since $(s,\half t)\in B(\delta_2)$,
applying (\ref{startit}) twice
we see that
\begin{eqnarray*}
g(\eta(s,\quart t)) &= &
\half g(\eta(s,\half t))
-\half g(\kappa(s,\half t))
-\half R(\omega(s, \half t))\\
&=&
\quart g(\eta(s, t)) -\quart g(\kappa(s,t))
-\quart R(\omega(s,t))
- \half g(\kappa(s,\half t))-\half R(\omega(s,\half t))
\,.
\end{eqnarray*}
Similarly, by a simple induction
\begin{equation}\label{beforedivis}
g(\eta(s,2^{-n}t)) = 2^{-n}g(\eta(s,t))-
\sum_{k=0}^{n-1} 2^{k-n}\Big(
g( \kappa(s,2^{-k}t))
+R(\omega(s,2^{-k}t))\Big)
\end{equation}
for all $(s,t)\in B(\delta_2)$ and $n\in \N$. Thus,
for all $(s,t)\in B(\delta_2)$ and $n\in \N$:
\begin{equation}\label{afterdivis}
\frac{g(\eta(s,2^{-n}t))}{2^{-n}}
= g(\eta(s,t))-
\sum_{k=0}^{n-1} 2^k \Big(
g( \kappa(s,2^{-k}t))
+R(\omega(s,2^{-k}t))\Big)\,.
\end{equation}
\end{numba}
Choose $\delta_3\in \;]0,\delta_2[$
such that $\delta_3\leq 1$;
then the closure $K:=\wb{B(\delta_3)}\sub B(\delta_2)$
is compact.
\begin{la}\label{lakapp}
There exists an absolutely convex,
bounded subset $B_1\not=\emptyset$ of $L(G)$
such that $\kappa(K)\sub L(G)_{B_1}$,
$\kappa|_K\!: K\to L(G)_{B_1}$
is Lipschitz continuous,
and $\|\kappa(s,t)\|_{B_1}\leq t^2$
for all $(s,t)\in K$.
In particular, $K_1:=\kappa(K)$ is compact
in $L(G)_{B_1}$.
\end{la}
\begin{proof}
Let $\rho_\kappa\!: B(\delta_1)^2\to L(G)$
be the first order Taylor remainder of~$\kappa$.
By Lemma~\ref{taylorR}, there exists an absolutely convex,
bounded subset $B_1\not=\emptyset$ of $L(G)$ such that
$\kappa(K)\sub L(G)_{B_1}$,
the map $\kappa|_K\!: K\to L(G)_{B_1}$
is Lipschitz continuous,
$\rho_\kappa(K\times K)\sub L(G)_{B_1}$,
and $\|\rho_\kappa(x,y)\|_{B_1}\leq (\|x-y\|_\infty)^2$
for all $x,y\in K$
(where $\|.\|_\infty$
is the maximum norm).
Note that $\frac{\partial \kappa}{\partial s}(s,0)=0$
for all $s\in \;]{-\delta_1},\delta_1[^{\ell-1}$
because $\kappa(\sbull,0)=0$.
Furthermore,
$\frac{\partial \kappa}{\partial t}(s,0)
=-2\cdot \half \cdot \frac{\partial \eta}{\partial t}(s,0)+
\frac{\partial \eta}{\partial t}(s,0)=0$
by the Chain Rule
(cf.\ (\ref{tauR})).
Hence $\kappa'(s,0)=0$ for all~$s$,
whence the Taylor expansion around $(s,0)$ yields
$\kappa(s,t)=\kappa(s,t)-\kappa(s,0)=\rho_\kappa((s,t),(s,0))$
for all $(s,t)\in K$. Therefore
$\|\kappa(s,t)\|_{B_1}=\|\rho_\kappa((s,t),(s,0))\|_{B_1}
\leq \|(s,t)-(s,0)\|_\infty^2=t^2$.
\end{proof}
\begin{numba}\label{estig}
By Lemmas~\ref{onbcompact} and~\ref{hoeloncp},
there is an absolutely
convex, bounded subset $B_2\not=\emptyset$
of $L(H)$ such that $g(K_1)\sub L(H)_{B_2}$
and $\|g(y)-g(x)\|_{B_2}\leq (\|y-x\|_{B_1})^\alpha$
for all $x,y\in K_1$.
\end{numba}
\begin{numba}\label{estiomeg}
The map $\omega$ being $h_\alpha$ and $K$ compact,
there exists an absolutely convex, bounded
subset $D \sub L(H)\times L(H)$ such that
$\omega(K)\sub (L(H)^2)_D$ and $\|\omega(y)-\omega(x)\|_D\leq
(\|y-x\|_\infty)^\alpha$ for all $x,y\in K$.
Then $K_2:= \omega(K)\sub S$ is compact in
$(L(H)^2)_D$.
\end{numba}
\begin{numba}\label{estiR}
Now Lemma~\ref{taylorR}
provides an absolutely convex, bounded subset
$B_3\not=\emptyset$ of $L(H)$
such that $R(K_2)\sub L(H)_{B_3}$,
$R|_{K_2}\!: K_2\to L(H)_{B_3}$ is Lipschitz continuous,
and $\|R(x)\|_{B_3}\leq (\|x\|_D)^2$ for all $x\in K_2$.
After replacing $B_2$ and $B_3$ with
an absolutely convex, bounded superset $B\sub L(H)$
(e.g., $B=\absconv(B_2\cup B_3)$), we may assume that
$B_2=B_3=B$.
\end{numba}
\begin{numba}\label{estR2}
Let $\rho_\eta\!: B(\delta)\times B(\delta)\to L(G)$
be the first order Taylor remainder of~$\eta$.
As $\eta$ is~$c^\infty$,
Lemma~\ref{taylorR}
provides
an absolutely convex, bounded subset
$B_4\not=\emptyset$ of $L(G)$
such that $\eta(K)\sub L(G)_{B_4}$,
$\eta|_K\!: K\to L(G)_{B_4}$ is Lipschitz continuous,
$\rho_\eta(K\times K)\sub
L(G)_{B_4}$,
and
\begin{equation}\label{thetaRem}
\|\rho_\eta(x,y)\|_{B_4}\leq (\|y-x\|_\infty)^2\quad
\mbox{for all $x,y\in K$.}
\end{equation}
Then $g$ being $h_\alpha$
and $\eta(K)\sub L(G)_{B_4}\cap U$
being compact,
after increasing~$B$
we may assume that
$g(\eta(K))\sub B\sub L(H)_B$ and
\begin{equation}\label{gontheta}
\|g(y)-g(x)\|_B\leq (\|y-x\|_{B_4})^\alpha
\quad\mbox{for all $x,y\in \eta(K)$}
\end{equation}
(Lemmas \ref{onbcompact} and \ref{hoeloncp}).
Then $g\circ \eta|_K\!: K\to L(H)_B$
is $H_\alpha$ (Lemma~\ref{laobvio}\,(b)).
After passing
to the closure, w.l.o.g.
$B$ is closed
in~$L(H)$ and so $(L(H)_B,\|.\|_B)$
is a Banach space ({\bf \ref{defnEB}}).
\end{numba}
We deduce estimates on the summands in (\ref{afterdivis})
(or multiples thereof) now.
\begin{numba}\label{estt1}
We have
$\|g(\kappa(s,2^{-k}t))\|_B
= \|g(\kappa(s,2^{-k}t))-g(\kappa(s,0))\|_B
\leq
\|\kappa(s,2^{-k}t)-\kappa(s,0)\|_{B_1}^\alpha
=\|\kappa(s,2^{-k}t)\|_{B_1}^\alpha\leq
2^{-2\alpha k}|t|^{2\alpha}$
for any $(s,t)\in K$ and $k\in \N_0$,
by {\bf \ref{estig}\/} and Lemma~\ref{lakapp}.
\end{numba}
\begin{numba}\label{estt2}
For any $(s,t)\in K$ and $k\in \N_0$,
we have
$\|R(\omega(s,2^{-k}t))\|_B\leq \|\omega(s,2^{-k}t)\|_D^2
=\|\omega(s,2^{-k}t)-\omega(s,0)\|_D^2\leq
2^{-2\alpha k}|t|^{2\alpha}$,
by {\bf \ref{estiR}\/} and {\bf \ref{estiomeg}}.
\end{numba}
\begin{la}\label{improexp}
If $\alpha\in \;]0,\half]$,
then $f$ is also $h_\beta$ for $\beta:=\frac{3}{2}$.
\end{la}
\begin{proof}
Let $\ell=2$.
For $\theta\!: \R^2\to G$ as in
Lemma~\ref{hoeldatpoint},
there is $\delta\in \;]0,\infty]$
such that $\theta(B(\delta))\sub U_1$.
Then {\bf \ref{makestarsh}}--{\bf\ref{estt2}} apply
to $\eta:=\phi\circ \theta|_{B(\delta)}$.
Let $I:=\;]{-\delta_3/2},\delta_3/2[$.
Since $\psi\circ f\circ \theta|_{B(\delta)}=
g\circ \eta$,
we have $\psi(f(\theta(I^2)))\sub g(\eta(K))\sub L(H)_B$
and $\psi\circ f\circ \theta|_{I^2}=g\circ \eta|_{I^2}\!:
I^2\to L(H)_B$ is continuous
(see {\bf \ref{estR2}}).
Thus conditions (a) and (b) of Lemma~\ref{hoeldatpoint}
are satisfied. To verify condition~(c),
let
$0\not=t\in I$
and $s\in I$.
There is $n\in \N$ with
$\half \delta_3 < |2^nt|\leq \delta_3$.
Let $t^*:=2^nt$; then $(s,t^*)\in K$ and thus
\begin{eqnarray*}
\|g(\eta(s,t)\|_B & = & \|g(\eta(s,2^{-n}t^*))\|_B\\
&\leq & 
2^{-n}\|g(\eta(s,t^*))\|_B+
\sum_{k=0}^{n-1} 2^{k-n}\Big(
\|g( \kappa(s,2^{-k}t^*))\|_B
+\|R(\omega(s,2^{-k}t^*))\|_B\Big)\\
&\leq &
2^{-n}+2 \sum_{k=0}^{n-1}
2^{k-n}2^{-2\alpha k}|t^*|^{2\alpha}
=
2^{-n}+2|t^*|^{2\alpha} 2^{-n} \sum_{k=0}^{n-1}
2^{(1-2 \alpha) k}\\
&\leq&
2^{-n}+2|t^*|^{2\alpha} 2^{-n} \sum_{k=0}^{n-1}
2^{(1-\frac{3}{2}\alpha) k}=
2^{-n}+\frac{2|t^*|^{2\alpha}2^{-n}}{2^{1-\frac{3}{2}\alpha}-1}
(2^{(1-\frac{3}{2} \alpha) n}-1)\\
&\leq &
2^{-\frac{3}{2}\alpha n}
+\frac{2|t^*|^{2\alpha}2^{-n}}{2^{1-\frac{3}{2}\alpha}-1}
2^{(1-\frac{3}{2} \alpha) n}
\leq c_1 (2^{-n})^{\frac{3}{2}\alpha}
<c_1(2/\delta_3)^{\frac{3}{2}\alpha}
|t|^{\frac{3}{2}\alpha}=c_2|t|^{\frac{3}{2}\alpha}
\end{eqnarray*}
with $c_1:=1+
\frac{2(\delta_3)^{2\alpha}}{2^{1-\frac{3}{2}\alpha}-1}$
and $c_2:=c_1(2/\delta_3)^{\frac{3}{2}\alpha}$.
Here, we used (\ref{beforedivis})
to pass to the second line,
then $g(\eta(K))\sub B$ (see {\bf \ref{estR2}}\/),
{\bf \ref{estt1}\/}
and {\bf \ref{estt2}\/}
to pass to the third.
We deduce that\linebreak
$\|\psi(f(\theta(s,t)))-\psi(f(\theta(s,0)))\|_B=\|\psi(f(\theta(s,t)))\|_B
=\|g(\eta(s,t))\|_B
\leq c_2 |t|^{\frac{3}{2}\alpha}$
for all $(s,t)\in I^2$
(for $t=0$, this is trivial).
Thus all conditions of Lemma~\ref{hoeldatpoint}
hold, and thus $f$ is~$h_{\frac{3}{2}\alpha}$.
\end{proof}
After replacing
$\alpha$ with $\beta:=\frac{3}{2}\alpha$
and repeating this process if necessary,
we see that $f$ is $h_\alpha$
with $\alpha\in \;]\half,1]$.
Thus, we may assume throughout the following
that $\alpha\in \;]\half,1]$.
We retain the setting and notations from
{\bf \ref{makestarsh}}--{\bf\ref{estt2}}.
\begin{la}\label{thelimit}
For each $(s,t)\in K$
such that $t\not=0$,
the limit
\begin{equation}\label{deflambd}
\lambda_\eta(s,t):=\lim_{n\to\infty}
\frac{g(\eta(s,2^{-n}t))}{2^{-n}t}
\end{equation}
exists in $(L(H)_B,\|.\|_B)$.
The convergence of
$\frac{g(\eta(s,2^{-n}t))}{2^{-n}t}$
in $L(H)_B$ is uniform
for $(s,t)$ in
$K^\times :=\{(s,t)\in K\!: t\not=0\}$.
The map $\lambda_\eta\!: K^\times \to L(H)_B$
is continuous, and
\begin{equation}\label{preindep}
\lambda_\eta(s,t)=\lambda_\eta(s,2^{-m}t)\quad\mbox{for all
$(s,t)\in K^\times$
and $m\in \Z$
such that $(s,2^{-m} t)\in K^\times$.}
\end{equation}
\end{la}
\begin{proof}
Given $\ve>0$,
choose
$N\in \N$ such that $c\,2^{-(2\alpha-1) N}\leq \ve$,
where $c:=\frac{2 (\delta_3)^{2\alpha-1}}{1-2^{-(2\alpha-1)}}$.
Using
(\ref{afterdivis}),
{\bf \ref{estt1}},
{\bf \ref{estt2}\/}
and the summation formula for the geometric
series, we obtain
for all $(s,t)\in K^\times$ and all $n,m\in \N$
such that $m\geq n\geq N$ the following estimates:
\begin{eqnarray*}
\left\| \frac{g(\eta(s,2^{-m}t))}{2^{-m}t}
- \frac{g(\eta(s,2^{-n}t))}{2^{-n}t}\right\|_B
\!\!\!\!\!\!& \leq &\!\! |t|^{-1}\! \sum_{k=n}^{m-1} 2^k \Big(
\|g(\kappa(s,2^{-k}t))\|_B+\|R(\omega(s,2^{-k}t))\|_B\Big)\\
\!\!\!\!\!\!&\leq & \!\! |t|^{-1}\! \sum_{k=n}^{m-1} 2\cdot 2^{-k}
\cdot 2^{-2\alpha k}|t|^{2\alpha}
= 2 |t|^{2\alpha-1}\sum_{k=n}^{m-1}2^{-(2\alpha-1) k}\\
\!\!\!\!\!\!&\leq &\!\! 2 (\delta_3)^{2\alpha-1}\!
\sum_{k=n}^{m-1} 2^{-(2\alpha-1) k}
\leq
c\,2^{-(2\alpha-1) n}\leq c\, 2^{-(2\alpha-1) N}\leq \ve.
\end{eqnarray*}
Thus $\big(\frac{g(\eta(s,2^{-n}t))}{2^{-n}t}\big)_{n\in \N}$
is a Cauchy sequence in the Banach space $(L(H)_B,\|.\|_B)$,
and hence the limit in (\ref{deflambd}) exists.
Letting $m\to \infty$ in the preceding inequalities,
we see that\linebreak
$\|\frac{g(\eta(s,2^{-n}t))}{2^{-n}t}-\lambda_\eta(s,t)\|_B
\leq \ve$
for all $(s,t)\in K^\times$ and $n\geq N$.
Hence $\mu_n\to \lambda_\eta$ uniformly,
where $\mu_n\!: K^\times \to L(H)_B$,
$\mu_n(s,t):=\frac{g(\eta(s,2^{-n}t))}{2^{-n}t}$.
Each $\mu_n$ being continuous
as a map into $L(H)_B$ (cf.\ {\bf \ref{estR2}}),
so is the uniform limit~$\lambda_\eta$.
To prove the final assertion,
we may assume that $m\geq 0$ (otherwise,
interchange the roles of $t$ and $2^{-m} t$).
Then $\lambda_\eta(s,t)
=\lim_{n\to \infty}
\frac{g(\eta(s,2^{-n}t))}{2^{-n}t}
=\lim_{n\to \infty}
\frac{g(\eta(s,2^{-(n+m)}t))}{2^{-(n+m)}t}
=\lim_{n\to \infty}
\frac{g(\eta(s,2^{-n}2^{-m}t))}{2^{-n}2^{-m}t}
=\lambda_\eta(s,2^{-m}t)$.
\end{proof}
We now specialize to $\ell=1$
for the rest of the proof (with
the exception of Lemma~\ref{Aisbd}).
Thus $\gamma:=\eta\!: \;]{-\delta},\delta[\,\to U \sub L(G)$
is a $c^\infty$-curve, with $\gamma(0)=0$.
\begin{la}\label{gamprime}
If $\gamma'(0)=0$,
then $\lambda_\gamma(t)=0\,$
for all $t\in K^\times=[{-\delta_3},\delta_3]\setminus\{0\}$.
\end{la}
\begin{proof}
If $\gamma'(0)=0$, then $\gamma(t)=\rho_\gamma(0,t)$
for each $t\in K$
(where $\rho_\gamma$ is the first order Taylor remainder
of $\gamma$).
Hence $\|g(\gamma(t))\|_B\leq \|\gamma(t)\|_{B_4}^\alpha
=\|\rho_\gamma(0,t)\|_{B_4}^\alpha\leq |t|^{2\alpha}$,
using (\ref{thetaRem}) and (\ref{gontheta}).
Therefore $\|t^{-1}g(\gamma(t))\|_B\leq |t|^{2\alpha-1}\to 0$
as $t\to 0$ and thus $\lim_{t\to 0}\frac{g(\gamma(t))}{t}=0$
in $L(H)_B$, entailing that $\lambda_\gamma(t)=\lim_{n\to \infty}
\frac{g(\gamma(2^{-n}t))}{2^{-n}t}=0$
for each $t\in K^\times$.
\end{proof}
To emphasize their dependence on $\gamma$,
we now write
$\delta_\gamma$, $\delta_{3,\gamma}$, $\theta_\gamma$
and $K_\gamma^\times$
for $\delta$, $\delta_3$, $\theta$ and~$K^\times$.
\begin{la}\label{samederiv}
Assume that $\gamma_1\!: \R\supseteq B(\delta_{\gamma_1})\to U$
and $\gamma_2\!: \R\supseteq B(\delta_{\gamma_2})\to U$
are $c^\infty$-curves such that $\gamma_1(0)=\gamma_2(0)=0$
and
$\gamma_1'(0)=\gamma_2'(0)$. Then
$\lambda_{\gamma_1}(t)=\lambda_{\gamma_2}(t)$
for all $t\in K_{\gamma_1}^\times\cap K_{\gamma_2}^\times$.
\end{la}
\begin{proof}
There is
$\delta\in \;]0,\min\{\delta_{\gamma_1},\delta_{\gamma_2}\}]$
such that the $c^\infty$-curve
$\xi \!: B(\delta)\to G$, $\xi(t)
:=\theta_{\gamma_1}(t)\theta_{\gamma_2}(t)^{-1}$
has image in~$U_1$.
Then $\gamma:=\phi\circ \xi\!: B(\delta)\to U$
is a $c^\infty$-curve such that
$\gamma(0)=0$
and $\gamma'(0)=\gamma_1'(0)-\gamma_2'(0)=0$. Hence
$\lambda_\gamma(t)=0$
for all $t\in K_\gamma^\times$, by Lemma~\ref{gamprime}.
The group multiplication and inversion in~$H$
being $c^\infty$, the set
$M_1:=\{(x,y)\in V_1\times V_1\!: xy^{-1}\in V_1\}$
is open in $V_1\times V_1\sub H\times H$
(equipped with the $c^\infty$-topology).
Thus $M:=(\psi\times \psi)(M_1)$
is $c^\infty$-open in $L(H)\times L(H)$, and the map
$\mu\!: M\to V$, $\mu(x,y):=\psi(\psi^{-1}(x)\psi^{-1}(y)^{-1})$
is $c^\infty$. There is $\ve\in \;]0,\min\{\delta_{3,\gamma_1},
\delta_{3,\gamma_2},
\delta_{3,\gamma}\}]$ such that $(g\circ \gamma_1,g\circ
\gamma_2)(B(\ve))\sub M$.
Set $J :=\{2^{-n}\!: n\in \N\}\cup\{0\}$.
Fix $t\in B(\ve)$ such that $t\not=0$.
Then $h_\gamma \!: J \to L(H)$, $h_\gamma(r):=
\frac{g(\gamma(rt))}{t}$
and the analogous functions $h_{\gamma_1}$
and $h_{\gamma_2}\!: J\to L(H)$
are bornologically differentiable at~$0$,
with $h_\gamma '(0)=\lambda_\gamma(t)$,
$h_{\gamma_1}'(0)=\lambda_{\gamma_1}(t)$
and $h_{\gamma_1}'(0)=\lambda_{\gamma_1}(t)$
(see Lemma~\ref{thelimit}).
Since $h_\gamma=t^{-1}\mu\circ (t h_{\gamma_1},t h_{\gamma_2})$,
we deduce from Lemma~\ref{bornolgood}
that
\begin{eqnarray*}
0 &= & \lambda_\gamma(t)= h_\gamma '(0)=t^{-1} d\mu\big((0,0),\,
(th_{\gamma_1}'(0), t h_{\gamma_2}'(0))\big)
=d\mu \big((0,0),\, (\lambda_{\gamma_1}(t),\lambda_{\gamma_2}(t))\big)\\
&=& \lambda_{\gamma_1}(t)-\lambda_{\gamma_2}(t)
\end{eqnarray*}
and thus $\lambda_{\gamma_1}(t)=\lambda_{\gamma_2}(t)$.
Now let $t\in K^\times_{\gamma_1}\cap K^\times_{\gamma_2}$.
There is $n\in \N$ such that $|2^{-n}t|<\ve$.
Then $0\not= 2^{-n}t\in B(\ve)$ and hence
$\lambda_{\gamma_1}(t)=\lambda_{\gamma_1}(2^{-n}t)=\lambda_{\gamma_2}(2^{-n}t)
=\lambda_{\gamma_2}(t)$ by the special case just treated
and (\ref{preindep}).
\end{proof}
As before, $\gamma\!: \;]{-\delta_\gamma},\delta_\gamma[\to U \sub L(G)$ is any
$c^\infty$-curve with $\gamma(0)=0$.
\begin{la}\label{laindep}
$\lambda_\gamma(t_1)=\lambda_\gamma(t_2)$
holds, for any $\, t_1,t_2\in K^\times$.
\end{la}
\begin{proof}
By continuity of $\lambda_\gamma$ (Lemma~\ref{thelimit}),
it suffices to show that $\lambda_\gamma(t_1)=\lambda_\gamma(t_2)$
for all $t_1,t_2\in K^\times_\gamma \cap \Q$.
Given such $t_1,t_2$,
there exist $m_1,m_2\in \Z\setminus\{0\}$
such that $m_1t_1=m_2t_2$.
For $\delta\in \;]0,1]$ sufficiently small,
the $c^\infty$-curve $\xi\!: B(\delta)\to G$,
$\xi(t):=\theta_\gamma(t_1t)^{m_1}\theta_\gamma(t_2 t)^{-m_2}$
has image in~$U_1$; we define $\eta:=\phi\circ \xi\!: B(\delta)\to U$.
Then $\eta(0)=1$ and $\eta'(0)=
m_1t_1\gamma'(0)-m_2t_2\gamma'(0)=0$
and hence $\lambda_\eta(t)=0$
for all $t\in K^\times_\eta$,
by Lemma~\ref{gamprime}.
Define $M_1:=\{(x,y)\in V_1\times V_1 \!: x^{m_1}y^{-m_2}\in V_1\}$,
$M:=(\psi\times \psi)(M_1)$
and $\mu\!: M\to V$, $\mu(x,y):=\psi(\psi^{-1}(x)^{m_1}\psi^{-1}(y)^{-m_2})$.
Then
$(g(\gamma(t_1 t)),g(\gamma(t_2t)))\in M$
for all $t\in B(\delta)$
and
$g(\eta(t))=\mu(g(\gamma(t_1 t)),g(\gamma(t_2t)))$
for all~$t\in B(\delta)$.
Pick $m\in \N$
such that $2^{-m}\in K^\times_\eta$, and fix $t:=2^{-m}$.
We consider
the three mappings
$h_1,h_2,h_3\!: J\to L(H)$
on $J:= \{ 2^{-n} \!: n\in \N\}\cup\{0\}$
given by $h_1(r):= \frac{g(\gamma(rtt_1))}{tt_1}$,
$h_2(r):=\frac{g(\gamma(rtt_2))}{tt_2}$ and
$h_3(r):=\frac{g(\eta(rt))}{t}$,
respectively.
Since $h_3(r)=t^{-1} \mu(tt_1 h_1(r), tt_2h_2(r))$
for all $r\in J$,
$h_1'(0)=\lambda_\gamma(tt_1)=\lambda_\gamma(t_1)$,
$h_2'(0)=\lambda_\gamma(tt_2)=\lambda_\gamma(t_2)$
and $h_3'(0)=\lambda_\eta(t)=0$,
we
deduce as in the proof of Lemma~\ref{samederiv} that
$0=\lambda_\eta(t)=
t^{-1} d\mu\big((0,0),\, (tt_1 \lambda_\gamma(t_1),
tt_2\lambda_\gamma(t_2))\big)
=m_1t_1\lambda_\gamma(t_1)-m_2t_2\lambda_\gamma(t_2)
=m_1t_1\big(\lambda_\gamma(t_1)-\lambda_\gamma(t_2)\big)$.
Hence $\lambda_\gamma(t_1)=\lambda_\gamma(t_2)$ indeed.
\end{proof}
By Lemma~\ref{laindep},
$\Lambda(\gamma):=\lambda_\gamma(t_0)$ for $t_0\in K_\gamma^\times$
is well defined, independent of the choice of~$t_0$.
\begin{la}\label{finallywelld}
Let $\gamma_1\!: B(\delta_{\gamma_1})\to U\sub L(G)$
and $\gamma_2\!: B(\delta_{\gamma_2})\to U$
be $c^\infty$-curves such that $\gamma_1(0)=\gamma_2(0)$.
If $\gamma_1'(0)=\gamma_2'(0)$,
then $\Lambda(\gamma_1)=\Lambda(\gamma_2)$.
\end{la}
\begin{proof}
For small $t_0\not=0$,
we have $\Lambda(\gamma_1)=\lambda_{\gamma_1}(t_0)=\lambda_{\gamma_2}(t_0)
=\Lambda(\gamma_2)$, by Lemma~\ref{samederiv}.
\end{proof}
Given $v\in L(G)$,
we choose
a $c^\infty$-curve $\gamma \!:B(\delta)\to U$
for some $\delta \in \;]0,\infty]$ 
such that $\gamma(0)=0$ and $\gamma'(0)=v$, and set
$A.v:=\Lambda(\gamma)\in L(H)$.
Lemma~\ref{finallywelld}
implies that $A.v$ is well defined,
independent~of~$\gamma$.
\begin{la}\label{Alin}
The map $A\!: L(G)\to L(H)$ is linear.
\end{la}
\begin{proof}
{\em $A$ is additive.}
Given $v_1,v_2\in L(G)$,
we choose
$\delta \in \;]0,\infty]$
and
$c^\infty$-curves $\gamma_j \!: B(\delta)\to U$
for $j\in \{1,2\}$
such that $\gamma_j(0)=0$ and $\gamma'_j(0)=v_j$.
After shrinking $\delta$, we may assume that
$\xi\!: B(\delta)\to G$, $\xi(t):=\phi^{-1}(\gamma_1(t))\phi^{-1}(\gamma_2(t))$
has image in~$U_1$; we define $\gamma\!: B(\delta)\to U$,
$\gamma(t):=\phi(\xi(t))$.
Then $\gamma(0)=0$ and $\gamma'(0)=\gamma_1'(0)+\gamma_2'(0)=v_1+v_2$.
Define $M_1:=\{(x,y)\in V_1\times V_1 \!: xy\in V_1\}$,
$M:=(\psi\times \psi)(M_1)$
and $\mu\!: M\to V$, $\mu(x,y):=\psi(\psi^{-1}(x)\psi^{-1}(y))$.
Then $\,\im(g\circ \gamma_1,g\circ \gamma_2)\sub M$
and $g(\gamma(t))=\mu(g(\gamma_1(t)), g(\gamma_2(t)))$
for all $t\in B(\delta)$.
Let $t\in K_\gamma^\times \cap
K_{\gamma_1}^\times \cap K^\times_{\gamma_2}$.
As in the proof of Lemma~\ref{samederiv},
we conclude that $\lambda_\gamma(t)=d\mu \big(
(0,0),\, (\lambda_{\gamma_1}(t), \lambda_{\gamma_2}(t))\big)
=\lambda_{\gamma_1}(t)+\lambda_{\gamma_2}(t)$.
Thus $A(v_1+v_2)=\Lambda(\gamma)=\lambda_\gamma(t)=
\lambda_{\gamma_1}(t)+\lambda_{\gamma_2}(t)=
\Lambda(\gamma_1)+\Lambda(\gamma_2)=
A(v_1)+A(v_2)$.\\[3mm]
{\em Homogeneity.}
Let $v\in L(G)$ and $a \in \R$.
If $a=0$, then $A(av)=A(0)=0=aA(v)$,
using that $A(0)=0$ as $A$ is a homomorphism
of additive groups.
If $a\not=0$,
pick a $c^\infty$-curve $\gamma\!:\R\to U$
such that $\gamma(0)=0$ and $\gamma'(0)=v$.
Then $\gamma_a\!: \R\to U$, $\gamma_a(t):=\gamma(at)$
is a $c^\infty$-curve with
$\gamma_a(0)=0$ and $\gamma'_a(0)=av$.
There is $t\in K_{\gamma_a}^\times$
such that $at\in K_\gamma^\times$.
Then $aA(v)=$ $a\Lambda(\gamma)=a\lambda_\gamma(at)
={\displaystyle \lim_{n\to\infty}}\frac{g(\gamma(2^{-n}at))}{2^{-n}t}
={\displaystyle \lim_{n\to\infty}}\frac{g(\gamma_a(2^{-n}t))}{2^{-n}t}
=\lambda_{\gamma_a}(t)=\Lambda(\gamma_a)=A(av)$.
\end{proof}
\begin{la}\label{Aisbd}
The linear map $A\!: L(G)\to L(H)$ is bounded.
\end{la}
\begin{proof}
By Lemma~\ref{bdlin},
it suffices to show
that $A\circ \gamma$ is continuous at~$0$
for each $c^\infty$-curve\linebreak
$\gamma\!: \R\to L(G)$
such that $\gamma(0)=0$.
Given such~$\gamma$,
there is $\delta>0$ with
$\;]{-\delta},\delta[\cdot \gamma(]{-\delta},\delta[)
\sub U$.
We define $\eta\!:
\R^2 \supseteq B(\delta)\to U$,
$\eta(s,t):=t \gamma(s)$
and $\eta_s\!: \;]{-\delta},\delta[\,\to U$,
$\eta_s(t):=\eta(s,t)=t\gamma(s)$ for $s\in \;]{-\delta},\delta[$.
Fix $t_0\in \;]0, \delta_{3,\eta}[$.
Then $\eta_s'(0)=\gamma(s)$
for all $s\in [{-\delta_{3,\eta}},\delta_{3,\eta}]$,
and apparently we can choose $\delta_{3,\eta_s}=\delta_{3,\eta}$.
Thus $A(\gamma(s))=\Lambda(\eta_s)=\lambda_\eta(s,t_0)$,
which is continuous in~$s$ as a map into $L(H)_B$
(see Lemma~\ref{thelimit}) and hence also as a map
into~$L(H)$, as desired.
\end{proof}
\begin{la}\label{ess2}
For every $c^\infty$-curve
$\gamma\!: \R \to U$ such that $\gamma(0)=0$,
the composition $g\circ \gamma$
is bornologically differentiable at~$0$,
with $(g\circ \gamma)'(0)=A.\gamma'(0)$.
\end{la}
\begin{proof}
We re-use the notations from {\bf \ref{makestarsh}}--{\bf \ref{estt2}}
and Lemma~\ref{thelimit}
for $\eta:=\gamma$.
Since $g(\gamma(K))\sub L(H)_B$
and $A.\gamma'(0)=\Lambda(\gamma)\in L(H)_B$,
we only need to show that $\,\lim_{n\to\infty}\frac{g(\gamma(t_n))}{t_n}$
$=\Lambda(\gamma)$ in $L(H)_B$,
for each sequence $(t_n)_{n\in \N}$
in $K^\times$ such that $t_n\to 0$.
We may assume that $|t_n|\leq \half \delta_3$
for all~$n$. For each $n$, there exists
a unique $m_n\in \N$ such that $t_n^*:=2^{m_n}t_n$
has absolute value $\half \delta_3<|t_n^*|\leq \delta_3$.
Define $\mu_n\!: K^\times\to L(H)_B$,
$\mu_n(t):=\frac{g(\gamma(2^{-n}t))}{2^{-n}t}$
for $n\in \N$, as in the proof of Lemma~\ref{thelimit}.
Since $\mu_n\to \lambda_\gamma$
uniformly by Lemma~\ref{thelimit},
given $\ve>0$ there exists
$n_0 \in \N$ such that
$\|\mu_n(t)-\lambda_\gamma(t)\|_B\leq \ve$
for all $n\geq n_0$ and
$t\in K^\times$.
Since $m_n\to\infty$ as $n\to\infty$,
there exists $n_1\in \N$ such that $m_n\geq n_0$
for all $n\geq n_1$.
For each $n\geq n_1$, we then have
$\big\|\frac{g(\gamma(t_n))}{t_n}-\Lambda(\gamma)\big\|_B
=\big\|\frac{g(\gamma(2^{-m_n}t_n^*))}{2^{-m_n}t_n^*}-\lambda_\gamma(t_n^*)\big\|_B
=\|\mu_{m_n}(t_n^*)-\lambda_\gamma(t_n^*)\|_B\leq \ve$.
Thus $\lim_{n\to\infty}\frac{g(\gamma(t_n))}{t_n}=\Lambda(\gamma)$
in $L(H)_B$ indeed.
\end{proof}
By Lemma~\ref{Aisbd} and Lemma~\ref{ess2},
$g=\psi\circ f\circ \phi^{-1}$
is bornologically curve
differentiable at~$0$.
Therefore $g$ is
curve differentiable at~$0$ (Lemma~\ref{strongstronger}),
and hence
$f$ is curve differentiable
(and bornologically curve differentiable)
at~$1$. As a consequence,
$f$ is $c^\infty$ (Lemma~\ref{suffs}).
This completes the proof of Theorem~\ref{main}.
\end{proof}
\appendix
\section{Proofs for Lemma~1.7 and Lemma~2.9}\label{appB}
\noindent
{\bf Proof of Lemma~\ref{BGNcalc}.}
Since $\mu\!: U\times E\times \R\to E$,
$\mu(x,y,t):=x+ty$ is a $c^\infty$-map
and hence continuous with respect to the $c^\infty$-topologies,
the preimage
$U^{[1]}=\mu^{-1}(U)$ is $c^\infty$-open in $U\times E\times \R$.
Let $\gamma=(\gamma_1,\gamma_2,\tau)\!: \R\to U^{[1]}$
be a smooth curve, with coordinates $\gamma_1\!:\R\to U$,
$\gamma_2\!: \R\to E$ and $\tau\!: \R\to \R$,
respectively.
Then
$I:=\tau^{-1}(\R\setminus\{0\})$
is an open subset of~$\R$.
To see that $f^{[1]}$ is $c^\infty$,
we have to show that $f^{[1]}\circ \gamma$
is smooth.
Clearly
$f^{[1]}\circ \gamma|_I$
is smooth, because this function
is composed of $c^\infty$-maps:
$\,f^{[1]}(\gamma(t))=\tau(t)^{-1}\big(f(\gamma_1(t)+\tau(t)\gamma_2(t))-
f(\gamma_1(t))\big)$
for $t\in I$.
Now assume that $t_0\in \R\setminus I$;
thus $\tau(t_0)=0$.
The map $h\!: \R^2\to E$, $h(t,s):=
\gamma_1(t)+s\tau(t)\gamma_2(t)$
being $c^\infty$, with $h(t_0,s)=\gamma_1(t_0)\in U$ for all~$s$,
we see that $h^{-1}(U)$ is an open neighbourhood
of $\{t_0\}\times [{-1},2]$ in~$\R^2$.
We therefore find an open neighbourhood
$J\sub \R$ of~$t_0$ such that $J\times [{-1},2]\sub
h^{-1}(U)$. Then $J\times \,]{-1},2[\,\to F$,
$(t,s)\mto
df\big(\gamma_1(t)+s\tau(t)\gamma_2(t), \, \gamma_2(t)\big)$
is smooth, and we have
\[
f^{[1]}(\gamma(t))=\int_0^1
df\Big(\gamma_1(t)+s\, \tau(t)\gamma_2(t), \; \gamma_2(t)\Big)\, ds
\quad \mbox{for $t\in J$.}
\]
Indeed, this formula
is obvious if $\tau(t)=0$;
if $\tau(t)\not=0$, it follows from the fundamental theorem
of calculus \cite[Cor.\,2.6\,(6)]{KaM}.
Being given by
a parameter-dependent integral with smooth integrand,
$f^{[1]}\circ \gamma|_J\!: J\to F$ is smooth
(cf.\ \cite[Prop.\,3.15]{KaM} or \cite[La.\,7.5]{BGN}).\vspace{3mm}\Punkt

\noindent
To facilitate a proof
of Lemma~\ref{onbcompact},
we first need to establish
a variant
of \cite[Prop.\,4.3.3]{FaK}:
\begin{la}\label{smoothsuff}
Let $E$ and $F$ be convenient vector spaces,
$U\sub E$ be $c^\infty$-open, and
$f\!: U\to F$ be~$h_\alpha$.
Then $f\circ \gamma\!: \R\to F$ is
$H_\alpha$, for every Lipschitz continuous curve $\gamma\!:
\R\to U$.
\end{la}
\begin{proof}
In view of Lemma~\ref{hoelvsbd}\,(b),
we may assume that $F=\R$.
The proof is by contraposition.
Thus, assume that $f\circ \eta$ is not $H_\alpha$
for some Lipschitz continuous curve $\eta\!: \R\to U$.
Then there exists $t_0\in \R$ such that $\Big\{
\frac{f(\eta(s))-f(\eta(t))}{|s-t|^\alpha}
\!: s,t\in I,\,s\not=t\Big\}$ is unbounded
for any neighbourhood $I$ of~$t_0$.
After translations, without loss of generality
$t_0=0$ and $\eta(t_0)=0$.
For $n\in \N$, we choose
$s_n,t_n \in [{-2^{-2n}},2^{-2n}]$
such that $t_n \not=s_n$ and $|s_n-t_n|^{-\alpha}
|f(\eta(s_n))-f(\eta(t_n))|\geq n2^{\alpha n}$.
We abbreviate $\sigma_n:=2^ns_n$, $\tau_n:=2^n t_n$
and define $\eta_n\!: \R\to E$,
$\eta_n(t):=\eta(t_n)+(t-t_n)\frac{\eta(s_n)-\eta(t_n)}{s_n-t_n}$
and $\gamma_n\!: \R\to E$,
$\gamma_n(t):=\eta_n(2^{-n}t)$.
Then
\[
\frac{|f(\gamma_n(\sigma_n))-f(\gamma_n(\tau_n))|}{|\sigma_n-\tau_n|^\alpha}
=\frac{|f(\eta(s_n))-f(\eta(t_n))|}{2^{n\alpha}|s_n-t_n|^\alpha}\geq n\,.
\]
Furthermore, $|\sigma_n|=2^n|s_n|\leq 2^{-n}$
and likewise $|\tau_n|\leq 2^{-n}$.
We claim that the sequence $(\gamma_n)_{n\in \N}$
in $C^\infty(\R,E)$ is fast
falling, in the sense of \cite[Defn.\,4.2.14]{FaK}.
By \cite[Prop.\,4.2.16]{FaK} (or \cite[Cor.\,12.3]{KaM}),
we only need to show that $(\lambda(\gamma_n(t)))_{n\in \N}$
is fast falling in~$\R$, for each $t\in [{-1},1]$
and each bounded linear functional $\lambda$ on~$E$.
Since $\eta$ is Lipschitz continuous,
so is $\lambda\circ \eta\!: \R\to\R$.
Using Lemma~\ref{hoeloncp}\,(a), we therefore find
$K\in [0,\infty[$ such that $|\lambda(\eta(r))-\lambda(\eta(s))|
\leq K\, |r-s|$ for all $r,s\in [{-1},1]$.
For $t\in [{-1},1]$, we obtain
\begin{eqnarray*}
|\lambda(\gamma_n(t))| &=& |\lambda(\eta_n(2^{-n}t))|
=\left| \lambda(\eta(t_n))+(2^{-n}t-t_n)
\frac{\lambda(\eta(s_n)-\eta(t_n))}{s_n-t_n}\right|\\
&\leq &
|\lambda(\eta(t_n))-\lambda(\eta(0))|
+|2^{-n}t-t_n|\cdot
\left|\frac{\lambda(\eta(s_n)-\eta(t_n))}{s_n-t_n}\right|\\
&\leq & K\, |t_n| + (2^{-n}+|t_n|)\, K
\leq (2^{-2n}+2^{-n}+2^{-2n})\, K\,,
\end{eqnarray*}
which is fast falling in~$\R$ as $n\to\infty$
(passing to the second line,
we used that $\lambda(\eta(0))=0$).
Hence indeed $(\gamma_n)_{n\in \N}$ is fast falling
in $C^\infty(\R,E)$. Applying the General Curve Lemma
\cite[Prop.\,4.2.15]{FaK} (or \cite[12.2]{KaM})
with $\ve_n:=2^{-n}$,
we get a smooth curve $\gamma\!: \R\to E$
and a convergent sequence $(r_n)_{n\in \N}$ of reals,
with limit $r:=\lim_{n\to\infty}r_n$,
such that $\gamma(r)=0$
and
$\gamma(r_n+t)=\gamma_n(t)$
for all $n\in \N$ and $t\in \R$ such that $|t|\leq \ve_n$.
Since $\gamma(r)=0=\eta(0)\in U$, the set
$J:=\gamma^{-1}(U)$ is an open neighbourhood
of~$r$, and thus $\gamma|_J\!: J\to U$
is a smooth curve in~$U$.
There is $N\in \N$ such that
$r_n+\tau_n, r_n+\sigma_n\in J$ for all $n\geq N$.
For any such $n$,
\[
\frac{|f(\gamma(r_n+\sigma_n)
-f(\gamma(r_n+\tau_n))|}{|\sigma_n-\tau_n|^\alpha}
=2^{-\alpha n}
\frac{|f(\eta(s_n))-f(\eta(t_n))|}{|s_n-t_n|^\alpha}
\geq n\, ,
\]
entailing that $\Big\{
\frac{f(\gamma(s))-f(\gamma(t))}{|s-t|^\alpha} \! :
s,t \in W, \, s\not=t \Big\} \sub \R$
is unbounded for each neighbourhood
$W\sub J$ of~$r$.
Hence $f\circ \gamma|_J$ is not $H_\alpha$
and hence $f$ is not $h_\alpha$
(cf.\ Remark~\ref{whoeldhoeld}\,(b)).
\end{proof}
{\bf Proof of Lemma~\ref{onbcompact}.}
The implication ``(b)$\impl$(a)''
can be proved like Lemma~\ref{comphoelder}
(and we shall not use it).
``(a)$\impl$(b)'':
The proof is by contraposition;
we assume that (b) is false and so
$f|_K\!: K\to F$ is not $H_\alpha$
for some $K$. Then $\lambda\circ f|_K$
is not $H_\alpha$ for some continuous
linear functional $\lambda\!: F\to \R$
(cf.\ Lemma~\ref{hoelvsbd}\,(b)).
If we can show that $\lambda\circ f$
is not $h_\alpha$, then neither is~$f$
(Remark~\ref{whoeldhoeld}\,(a)).
Hence $F=\R$ without loss of generality.
As we assume that $f|_K$ is not $H_\alpha$,
for each $n\in \N$
we find elements $x_n,y_n\in K$
such that $\|y_n-x_n\|_B < 1/n^2$
and
$|f(y_n)-f(x_n)|\geq n(\|y_n-x_n\|_B)^\alpha$.
Using that $K$ is compact and metrizable,
after passing to subsequences
we may assume that both $x_n$ and $y_n$
converge to some $x\in K$,
and $\|y_n-x\|_B, \|x_n-x\|_B < 1/n^2$.
We now consider the curve
$\gamma\!: \R\to E_B$
defined as follows:
$\gamma(t):=x_1$ if $t\leq 0$;
$\gamma$~runs with constant velocity $\frac{y_1-x_1}{\|y_1-x_1\|_B}$
from $x_1$ to $y_1$ if $t_1:=0\leq t\leq \|y_1-x_1\|_B=:s_1$;
$\gamma$~runs with constant velocity
$\frac{x_2-y_1}{\|x_2-y_1\|_B}$
from $y_1$ to $x_2$
if $s_1\leq t\leq s_1+\|x_2-y_1\|_B=: t_2$,
and so on. Since
$t_\infty:=\sum_{n=1}^\infty
\|y_n-x_n\|_B +\sum_{n=1}^\infty \|x_{n+1}-y_n\|_B$
is finite, $\gamma(t)$ tends to $x$
as $t$ increases towards $t_\infty$;
so we define $\gamma(t):=x$ for $t\geq t_\infty$.
By construction, we have
$\|\gamma(s)-\gamma(t)\|_B\leq |s-t|$
for all $s,t\in \R$,
and thus $\gamma$ is Lipschitz continuous.
Since $\gamma(t_\infty)\in K\sub U$, the map
$\gamma\!: \R\to E_B$
is continuous, and $U\cap E_B$ is open in~$E_B$,
we deduce that
$J:=\gamma^{-1}(U)$ is an open neighbourhood of~$t_\infty$
in~$\R$. There is $N\in \N$ such that
$s_n, t_n\in J$ for all $n\geq N$.
For any such~$n$,
we have $\|\gamma(s_n)-\gamma(t_n)\|_B=\|y_n-x_n\|_B=|s_n-t_n|$
and hence $|f(\gamma(s_n))-f(\gamma(t_n))|
=|f(y_n)-f(x_n)| \geq n(\|y_n-x_n\|_B)^\alpha=
n|s_n-t_n|^\alpha$.
Since $t_n,s_n\to t_\infty$,
this implies that the set
$\Big\{\frac{f(\gamma(s))-f(\gamma(t))}{|s-t|^\alpha}\!: s,t\in W, s\not=t\Big\}\sub \R$
is unbounded for each neighbourhood
$W\sub J$ of $t_\infty$,
whence $f\circ \gamma|_J$ is not $H_\alpha$.
Therefore, by Lemma~\ref{smoothsuff},
$f$ is not $h_\alpha$.
The proof
is complete.\Punkt
\noindent
{\footnotesize
{\bf Helge Gl\"{o}ckner}, TU~Darmstadt, FB~Mathematik~AG~5,
Schlossgartenstr.\,7, 64289 Darmstadt, Germany.\\
E-Mail: gloeckner@mathematik.tu-darmstadt.de}

\begin{thebibliography}{10}
%
\bibitem{BGN} Bertram, W., H. Gl\"{o}ckner and K.-H. Neeb,
{\em Differential calculus over general base fields and rings},
to appear in Expo.\ Math.; cf.\ arXiv:math.GM/0303300\,.
%
%
\bibitem{Bil} Biller, H., {\em The exponential law
for smooth functions}, Manuscript, July 2002.
%
%
\bibitem{Bom} Boman, J., {\em Differentiability of a
function and of its compositions with functions of
one variable}, Math.\ Scand.\ {\bf 20}\,(1967), 249--268.
%
%
\bibitem{Bou} Bourbaki, N., ``Topological Vector Spaces,
Chapters~1--5,'' Springer-Verlag,
1987.
%
%
\bibitem{FaF} Faure, C.-A. and A. Fr\"{o}licher,
{\em H\"{o}lder differentiable maps and their
function spaces},
pp.\,135--142 in: ``Categorical Topology and its
Relation to Analysis, Algebra and Combinatorics''
(Prague, 1988), World Sci.\ Publ., Teaneck, NJ, 1989.
%
%
\bibitem{FaK} Fr\"{o}licher, A. and A. Kriegl,
``Linear Spaces and Differentiation Theory,''
Wiley-Interscience, Chichester, 1988.
%
%
\bibitem{RES} Gl\"{o}ckner, H.,
{\em Infinite-dimensional Lie groups without completeness
restrictions}, pp.\,43--59 in
Strasburger, A. et al.\ (Eds.),
``Geometry and Analysis
on Finite- and Infinite-Dimensional Lie Groups,''
Banach Center Publications~{\bf 55},
Warsaw, 2002.
%
%
\bibitem{HOL} -----, {\em H\"{o}lder continuous
homomorphisms between infinite-dimensional Lie groups are smooth},
TU Darmstadt Preprint {\bf 2327}, March 2004; also arXiv:math.GR/0403251\,.
%
%
\bibitem{KaM} Kriegl, A. and P.\,W. Michor,
``The Convenient Setting of Global Analysis,''
Math.\ Surveys and Monographs {\bf 53},
AMS, Providence, 1997.
%
%
\bibitem{Mi1} Milnor, J., {\em On infinite dimensional Lie groups},
Preprint, Institute for Advanced Study, Princeton, 1982.
%
%
\end{thebibliography}
\end{document}